\documentclass[final,3p,times]{elsarticle}

\usepackage[lined,ruled]{algorithm2e}
\usepackage{amsfonts}
\usepackage{amsmath}
\usepackage{amssymb}
\usepackage{caption}
\usepackage{color}
\usepackage{epsfig}
\usepackage{enumerate}
\usepackage{float}
\usepackage{graphicx}
\usepackage{subcaption}
\usepackage{lineno}
\usepackage{multirow}
\usepackage{natbib} 
\usepackage{placeins}
\usepackage{setspace}
\usepackage{stmaryrd}
\usepackage{url}
\usepackage{verbatim}

    \newcommand{\subf}[2]{%
      {\small\begin{tabular}[t]{@{}c@{}}
       \mbox{}\\[-\ht\strutbox]
       #1\\#2
       \end{tabular}}%
    }
    
    \newtheorem{theorem}{Theorem}[section]
    \newtheorem{corollary}{Corollary}
    \newdefinition{definition}{Definition}[section]
    
    \newcommand{\norm}[1]{\left\Vert#1\right\Vert}
    
    \newcommand{\set}[1]{\left\{#1\right\}}
    \newcommand{\mat}[2]{\left[\begin{array}{#1}#2\\ \end{array}\right]}
    
    \DeclareMathOperator*{\argmax}{arg\,max}

    \newcommand{\initialLabel}{\text{init}}
    \newcommand{\updatedLabel}{\text{update}}
    \newcommand{\observedLabel}{\text{obs}}
    \newcommand{\predictedLabel}{\text{pred}}
    
    \newcommand{\RR}{\mathbb{R}}
    \newcommand{\dspace}{{\mathcal{D}}}
    \newcommand{\pspace}{{\Lambda}}
    
    \newcommand{\param}{\lambda} 			
    \newcommand{\paramref}{\param^\dagger}          
    \newcommand{\mudpt}{\param^{\text{MUD}}}        
    \newcommand{\initialMean}{\param_0}             
    \newcommand{\observedMean}{\mathbf{y}}          
    
    \newcommand{\PP}{\mathbb{P}}
    \DeclareMathOperator*{\initialP}{\PP_\initialLabel}
    
    \DeclareMathOperator*{\observedP}{\PP_\observedLabel}
    \DeclareMathOperator*{\predictedP}{\PP_\predictedLabel}
    
    \newcommand{\pp}{\pi}
    \DeclareMathOperator*{\initial}{\pp_\initialLabel}
    \DeclareMathOperator*{\updated}{\pp_\updatedLabel}
    \DeclareMathOperator*{\observed}{\pp_\observedLabel}
    \DeclareMathOperator*{\predicted}{\pp_\predictedLabel}
    
    \newcommand{\cov}{\Sigma}
    \newcommand{\initialCov}{\cov_\initialLabel}
    \newcommand{\updatedCov}{\cov_\updatedLabel}
    \newcommand{\observedCov}{\cov_\observedLabel}
    \newcommand{\predictedCov}{\cov_\predictedLabel}

\journal{Computer Methods In Applied Mechanics and Engineering }

\begin{document}

\begin{frontmatter}

\title{Parameter Estimation with Maximal Updated Densities}

\author[mm]{Michael Pilosov}
\ead{michael.pilosov@ucdenver.edu}
\author[oden]{Carlos del-Castillo-Negrete}
\ead{carlos.del-castillo-negrete@utexas.edu}
\author[snl]{Tian Yu Yen}
\ead{tyen@sandia.gov}
\author[cud]{Troy Butler}
\ead{troy.butler@ucdenver.edu}
\author[oden]{Clint Dawson}
\ead{clint.dawson@utexas.edu}
\affiliation[mm]{organization={Mind the Math LLC},
            city={Denver},
            postcode={80203},
            state={CO},
            country={USA}}
\affiliation[oden]{organization={Oden Institute for Computational Engineering and Sciences, University of Texas at Austin},
            city={Austin},
            postcode={78712},
            state={TX},
            country={USA}}
\affiliation[snl]{organization={Sandia National Laboratories, Center for Computing Research},
            city={Albuquerque},
            postcode={87175},
            state={NM},
            country={USA}}
\affiliation[cud]{organization={Department of Mathematical and Statistical Sciences, University of Colorado Denver},
            city={Denver},
            postcode={80202},
            state={CO},
            country={USA}}


\begin{abstract}
    A recently developed measure-theoretic framework solves a stochastic inverse problem (SIP) for models where uncertainties in model output data are predominantly due to aleatoric (i.e., irreducible) uncertainties in model inputs (i.e., parameters).
    The subsequent inferential target is a distribution on parameters.
    Another type of inverse problem is to quantify uncertainties in estimates of ``true'' parameter values under the assumption that such uncertainties should be reduced as more data are incorporated into the problem, i.e., the uncertainty is considered epistemic.
    A major contribution of this work is the formulation and solution of such a parameter identification problem (PIP) within the measure-theoretic framework developed for the SIP.
    The approach is novel in that it utilizes a solution to a stochastic forward problem (SFP) to update an initial density only in the parameter directions informed by the model output data.
    In other words, this method performs ``selective regularization'' only in the parameter directions not informed by data.
    The solution is defined by a maximal updated density (MUD) point where the updated density defines the measure-theoretic solution to the PIP. 
    Another significant contribution of this work is the full theory of existence and uniqueness of MUD points for linear maps with Gaussian distributions.
    Data-constructed Quantity of Interest (QoI) maps are also presented and analyzed for solving the PIP within this measure-theoretic framework as a means of reducing uncertainties in the MUD estimate. 
    We conclude with a demonstration of the general applicability of the method on two problems involving either spatial or temporal data for estimating uncertain model parameters.
    The first problem utilizes spatial data from a stationary partial differential equation to produce a MUD estimate of an uncertain boundary condition.
    The second problem utilizes temporal data obtained from the state-of-the-art ADvanced CIRCulation (ADCIRC) model to obtain a MUD estimate of uncertain wind drag coefficients for a simulated extreme weather event near the Shinnecock Inlet located in the Outer Barrier of Long Island, NY, USA.
\end{abstract}





\end{frontmatter}



\section{Introduction}\label{sec:intro}

We increasingly rely on quantitative predictions from computational models of physical systems to inform engineering design, predict the behavior of physical systems, and even shape public policy, e.g., see \cite{omalley2014combined,omalley2015bayesianinformationgap, mattisParameterEstimationPrediction2015, harp2013contaminant}, for just a few such examples.
It is therefore more important than ever to quantify---and reduce whenever possible---the uncertainties impacting such models.

Unfortunately, many key characteristics governing system behavior, described as model inputs (referred to here as parameters), are often hidden from direct observation.
When observable model output data associated with some quantities of interest (QoI) are sensitive to variations in these parameters, we formulate and solve inverse problems using the output data for QoI to quantify uncertainties in parameters.
Inverse problems therefore play a vital role in the uncertainty quantification (UQ) community.

The type of inverse problem formulated and methods utilized for solving it fundamentally rely upon the type of uncertainties to be quantified.
In the UQ community, uncertainties are categorized as being either aleatoric (i.e., irreducible) or epistemic (i.e., reducible).
We refer to an inverse UQ problem involving aleatoric uncertainties as a stochastic inverse problem (SIP) whereas an inverse UQ problem involving epistemic uncertainties is referred to as a parameter identification problem (PIP). 
The solution to a SIP is a probability density (or, more generally, a probability measure) that characterizes the differences in relative likelihoods of parameters given a density characterizing the inherent variability of observable data.
However, the solution to the PIP is an estimate of a single ``true'' parameter, which may be obtained, for example, by maximizing a density that quantifies uncertainty in the estimate of this unknown, but fixed, parameter value due to finite observable data.

This brings us the contributions of this work, which we enumerate below for ease of reference.
\begin{enumerate}
	\item At a high-level, a major contribution of this work is the formulation and analysis of solutions to a PIP utilizing a novel measure-theoretic method originally developed for solving a SIP.
	\item In order to make direct comparisons to other formulations and solutions of the PIP based on popular Bayesian and least squares approaches, a full ``linear Gaussian'' theory of this new approach is developed. 
	\begin{itemize}
		\item As suggested by the subsequent comparison examples, this new approach is widely applicable to any PIP where Bayesian or least squares methods are regularly applied and has some potential advantages over solutions obtained via those methods.
		\item As the theory and examples demonstrate, the approach developed in this work applies ``selective regularization'' in the parameter space. Specifically, an initial characterization of parameter uncertainty does not impact the updated characterizations of uncertainty in the directions that the mapping deems ``informed'' by the data.
	\end{itemize}
	\item A theoretical analysis for data-constructed residual QoI maps within this framework is developed including a guarantee of existence, uniqueness, and convergence of solutions to the PIP as more data are collected for a fixed linear measurement.
	\item We conclude this work by presenting an alternative data-constructed residual QoI map for general (potentially nonlinear) spatial-temporal measurements based on a principal component analysis (PCA).
\end{enumerate}

Below, we summarize the organization of the paper as well as key results that are connected to the contributions enumerated above. \begin{description}
	\item[Section~\ref{sec:compare}:] This provides a brief literature review, motivation, and context for both the SIP and PIP, which sets the stage for contribution 1 of this work. 
		\begin{itemize}
			\item The formal definitions of the SIP and PIP are provided along with an illustrative example that compares the solutions to these problems.
			\item The example motivates the use of the {\em maximal updated density} (MUD) estimate in the measure-theoretic approach as an alternative to the Bayesian {\em maximum a posteriori} (MAP) point for solving the PIP. 
		\end{itemize}
	\item[Section~\ref{sec:data-maps}:] The first data-constructed Quantity of Interest (QoI) map considered in this work is defined as a weighted mean error map. This sets the stage for both contributions 2 and 3 of this work.
		\begin{itemize}
			\item The concept of statistical sufficiency motivates the construction of this map.
			\item This type of map is utilized for scenarios where multiple (i.e., repeated) noisy data are available for general types of measurements and motivates the development of the linear Gaussian theory (i.e., contribution 2).
		\end{itemize}
	\item[Section~\ref{sec:estimation}:] This is the heart of contributions 2 and 3 as we prove existence and uniqueness for MUD points as well as a convergence result for MUD estimates based on the first type of data-constructed QoI map (i.e., the weighted mean error map).
		\begin{itemize}
			\item It is shown that once a sufficient number of data are collected for each measurement, the existence and uniqueness of MUD points are guaranteed for weighted mean error maps.
			\item The variance in these estimates is shown to vanish in the limit of infinite data {\em only in the directions of the parameter space informed by the measurement operator}. In other words, selective regularization is applied in MUD estimation.
			\item A closed-form expression for the MUD point is provided along with conceptual and quantitative comparisons to MAP and least squares parameter estimates. The quantitative comparisons of these parameter estimates are performed for various randomly constructed linear maps on $100$-dimensional parameter spaces.
		\end{itemize}	 
	\item[Section~\ref{sec:PCA}:] This sets up contribution 4 of this work as we introduce an alternative data-constructed QoI based on a principal component analysis (PCA) that is easily applied to spatial and temporal data. Unlike the approach in Section~\ref{sec:data-maps}, this PCA based QoI approach does not require the availability of multiple (i.e., repeated) data for each measurement.
	\item[Section~\ref{sec:PCA_examples}:] This concludes contribution 4 of this work as we demonstrate the effectiveness of the PCA based QoI maps at producing accurate MUD estimates for two types of problems.
		\begin{itemize}
			\item The first example utilizes spatial data to form a PCA based QoI map to estimate an uncertain boundary condition for a two-dimensional stationary partial differential equation.
			\item The second example utilizes temporal data of free surface elevations to estimate wind drag coefficients for the state-of-the-art ADvanced CIRCulation (ADCIRC) model used to simulate an extreme weather event near the Shinnecock Inlet located in the Outer Barrier of Long Island, NY, USA. 
	\end{itemize}	
	\item[Section~\ref{sec:conclusions}:] Concluding remarks and future directions for this work are provided.
	\item[\ref{app:mud_deriv}:] The detailed derivation is provided for the closed-form expression of the MUD point for linear Gaussian maps that is utilized in Section~\ref{sec:estimation}. 
	\item[\ref{app:A}:] In the interest of scientific reproducibility, the details on obtaining the code and data utilized throughout this manuscript are provided. 
\end{description}




\section{Comparing Inverse Problems and Solutions}\label{sec:compare}

We begin with a formal discussion of the SIP and PIP in order to provide a brief literature review on these problems and their solutions.
This sets the stage for the remainder of this work and the contributions therein. 

The mapping from model parameters to observable model outputs defines what we refer to as a QoI map.
We assume the model parameters are hidden from direct observation and must be inferred from observable data associated with the QoI. 
When model parameters possess aleatoric uncertainties, e.g., due to naturally occurring variability in system inputs, then the specification of a probability measure quantifying uncertainties in the QoI data leads to the formulation of a SIP.
A solution to the SIP is given by a pullback of this probability measure onto the space of parameters.
We refer to a pullback measure as a data-consistent solution since its push-forward through the QoI map matches the probability measure quantifying uncertainties in the QoI data.

While it is possible to construct explicit approximations to data-consistent measures in terms of estimating measurable events and their probabilities in the parameter space (e.g., see \cite{butler2014measuretheoretic}), such ``set-based'' approximations become computationally intractable for high-dimensional parameter spaces or geometrically complex and/or computationally expensive QoI maps.
A recently developed density-based approach \citep{butlerCombiningPushForwardMeasures2018, butler2018convergence, butler2020dataconsistent} solves the SIP in a novel way by first solving a stochastic forward problem (SFP).
Specifically, an {\em initial} probability measure is first specified on the parameters to encode any prior knowledge of parameter variability.
Then, an SFP is solved where the push-forward of the initial probability measure is used to define a {\em predicted} probability measure on the QoI.
The discrepancy between the predicted and {\em observed} probability measures on the QoI, expressed as a ratio of probability density functions (more generally, Radon-Nikodym derivatives), is then used to {\em update} the initial probability density.
The {\em updated} probability measure associated with this density is then data-consistent.
Moreover, the updates to the initial probability measure only occur in directions informed by the QoI.
In other words, the initial probability measure serves to regularize the space of all pullback measures solving the SIP to produce a unique solution.

The SIP and its solution methodologies are based on rigorous measure theory using the Disintegration Theorem \citep{dellacherieProbabilitiesPotential1978, chang1997conditioning} as the central tool in establishing existence, uniqueness, and stability of solutions.
Updated probability measures often have complex structures that are not well approximated by a family of parametrically defined distributions (e.g., Gaussian).
This attribute of the solution further distinguishes this measure-theoretic approach from typical Bayesian-inspired approaches, e.g., Hierarchical Bayesian methods \citep{wang2004hierarchical, wikle1998hierarchical, Smith, tarantola2005inverse}, that specify prior distributions from a parametric family of distributions along with additional a priori assumed distributions on the hyper-parameters introduced by this parametric family (e.g., the means and variances of a Gaussian).
Subsequently, solutions to the SIP using Bayesian approaches will not, in general, produce solutions (defined as posterior distributions) whose push-forward matches the observed distribution.
In fact, the push-forward of the posterior is not even of general interest in most Bayesian paradigms.
Instead, the posterior predictive, which defines the distribution of possible unobserved values is of central interest \citep{Smith}.
The posterior predictive is constructed as a conditional distribution on the observations but makes practical use of the posterior through a marginalization.
These differences are not surprising when one considers that the Bayesian inverse problem that is perhaps most familiar in the UQ community solves an inverse problem involving epistemic uncertainty, as we describe and expand upon below.

In a typical Bayesian framework \citep{calvetti2014inverse, fitzpatrick1991bayesian, bui-thanh2014analysis, cockayneoatessullivangirolami, gelman2013bayesian, jaynes2003probability, kennedy2001bayesian, tarantola2005inverse, marzouk2007stochastic, cotter2010approximation, stark2010primer, ernst2014bayesian, cui2015scalable, Stuart10}, one of the initial assumptions is that data obtained on a QoI are polluted by measurement error, i.e., the data are ``noisy.''
Measurement errors can theoretically be reduced using improved measurement instruments (i.e., they are epistemic in nature).
A data-likelihood function is used to express the relative likelihoods that all of the data came from a particular choice of the parameter.
Encoding any initial assumptions about which parameters are more likely than others as a prior density allows the formal construction of a posterior density as a conditional density that describes the difference in relative likelihoods of any parameter value given the data.

It is common to use specific point estimators such as the maximum a posteriori (MAP) point given by the mode of the posterior as the actual solution to the inverse problem.
The posterior is then re-interpreted as providing descriptions of uncertainty in that specific point estimate.
The Bernstein-von Mises theorem \citep{vonmises1941discussion} provides conditions under which the posterior will become concentrated around the single true parameter in the limit of infinite data \citep{Smith}.
In other words, the use of such point estimates is useful when solving a PIP as the uncertainty in such estimates, usually quantified via the covariance of the posterior, ``shrinks'' as more data are incorporated. 

The Bayesian framework is fundamentally designed to solve a PIP while the measure-theoretic framework as presented in \cite{butlerCombiningPushForwardMeasures2018, butler2018convergence, butler2020dataconsistent} is instead designed to solve a SIP.
This difference in methods is not semantic as the SIP and PIP themselves arise under fundamentally different conditions. 
At a more conceptual level, a SIP arises more naturally in the context of applications where the objective is to quantify uncertainties across many instantiations of a system whereas a PIP arises more naturally in the context of applications where the objective is to quantify uncertainties in a specific realization of a system. 
For instance, in \cite{butler2020what}, SIPs and their data-consistent solutions are considered in the context of both a verification \& validation and quality control problem in the design and manufacturing of thin elastic membranes.
Separately, the work of \cite{tran2021solving} utilized data-consistent solutions to learn a distribution of microstructure parameters in the context of computational materials science.

To summarize at a high-level, the typical Bayesian approach to an inverse problem focuses on first modeling epistemic uncertainties in data on a QoI obtained from a true, but unknown, parameter value, which we denote by $\paramref$.
The objective is to then estimate this specific $\paramref$ and quantify uncertainties in this point estimate. 
This is in contrast to the SIP and its data-consistent solutions that are defined as pullback measures of an observed probability measure on the QoI.
To help build intuition about these differences, we summarize key technical details about these inverse problems and their solutions before presenting an example that highlights differences in solutions and the potential of using the data-consistent approach for solving the PIP.



\subsection{Terminology, notation, and the inverse problems}

To make comparisons more clear, we first introduce shared notation between the SIP and PIP.
Denote by $\pspace$ the space of (input) parameters for the model.
Denote by $Q$ the (potentially vector-valued) QoI map from the parameter space, $\pspace$, to the data space defined by $\dspace:= \set{Q(\lambda)\, : \, \lambda\in\pspace}$.
Note that while any given component of the QoI map is generally a composition of a functional applied to the solution of the model, we use the notation $Q(\lambda)$ to emphasize the dependence of the QoI on the choice of parameter.
For simplicity in presentation, we assume $\pspace \subseteq \mathbb{R}^p$ and $\mathcal{D} \subset \mathbb{R}^m$ for finite $p$ and $m$.
We use $Q^{-1}(E)$ for any $E\subset\mathcal{D}$ to denote the {\em preimage} of $E$, i.e., $Q^{-1}(E)=\set{\lambda\in\pspace \, : \, Q(\lambda)\in E}$.
Unless otherwise specified, we assume that $\pspace$ and $\dspace$ are equipped with (Borel) $\sigma$-algebras to define measurable spaces, $Q$ is a measurable map between these spaces, and that subsets of these spaces are taken from these $\sigma$-algebras. 

\subsubsection{The stochastic inverse problem (SIP) and parameter identification problem (PIP)}\label{subsec:defns}

We first define the concept of push-forward measures as solutions to the SFP mentioned above, which also helps frame the SIP more clearly as the direct inversion of the SFP.
\begin{definition}[Stochastic Forward Problem (SFP)]\label{def:forward-problem}
  Given an initial (i.e., initially specified) probability measure $\initialP$ on $\pspace$, the SFP is to determine the push-forward probability measure $\predictedP$ on $\dspace$, defined by
\begin{linenomath*}
\begin{equation*}
\predictedP(E) := \initialP(Q^{-1}(E)),
\end{equation*}
\end{linenomath*}
for all events $E\subset\dspace$.
\end{definition}

We often refer to the push-forward of the initial measure as the {\em predicted} measure since it may be constructed before any observed data are known.
This also helps to distinguish it from the {\em observed} measure used in the formulation of the SIP.
\begin{definition}[Stochastic Inverse Problem (SIP)]\label{def:inverse-problem}
Given an observed probability measure $\observedP$ on $\dspace$, the SIP is to determine a pullback probability measure, denoted $\mathbb{P}_\pspace$, on $\pspace$, which is data-consistent in the sense that
\begin{linenomath*}
\begin{equation}\label{eq:data-consistent}
\mathbb{P}_\pspace(Q^{-1}(E)) = \observedP(E),
\end{equation}
\end{linenomath*}
for all events $E\subset\dspace$.
\end{definition}

Unless the map $Q$ is a bijection, we do not expect that there is a unique $\mathbb{P}_\pspace$ solving the SIP, but rather there is a class of pullback measures that solve the SIP.
In \cite{butler2014measuretheoretic}, a disintegration theorem \cite{chang1997conditioning} along with an ansatz is used to establish the existence of solutions to the SIP that are unique up to the choice of ansatz.
An algorithm is provided in \cite{butler2014measuretheoretic} for explicitly approximating pullback measures by applying a specified ansatz to approximations of contour events, i.e., approximations of $Q^{-1}(E_i)$ where $\set{E_i}_{i\in\mathcal{I}}$ is a partition of $\dspace$.

In \cite{butlerCombiningPushForwardMeasures2018}, an alternative density-based approach is presented that is computationally simpler to implement, scales well with increasing parameter dimension, and is stable with respect to perturbations in the initial and observed probability measures.
We refer the interested reader to \cite{butlerCombiningPushForwardMeasures2018} for the theoretical and algorithmic details.
Here, we summarize the density-based solution to the SIP as:
\begin{linenomath*}
\begin{equation}\label{eq:updated-pdf}
	\updated(\param) := \initial(\param)\frac{\observed(Q(\param))}{\predicted(Q(\param))}.
\end{equation}
\end{linenomath*}
The densities (i.e., Radon-Nikodym derivatives) $\initial$ and $\observed$ are associated with the specification of $\initialP$ and $\observedP$, respectively.
The density $\predicted$ is associated with the predicted measure $\predictedP$.
The solution, $\updated$, is referred to as an updated density because it is defined as a multiplicative update to the initial density, $\initial$.
It is worth emphasizing that constructing $\updated$ as a solution to the SIP requires a solution to the SFP.
If $\initial$ and $\observed$ are prescribed, then the solution only requires the construction of $\predicted$, which is usually achieved via non-parametric means (see \cite{butlerCombiningPushForwardMeasures2018} for details). 
As a consequence, when $m<p$, the solution to the SIP is obtained by solving a forward UQ problem on a lower-dimensional space. 

In order to ensure that $\updated$ is in fact a density, a predictability assumption is required \cite{butlerCombiningPushForwardMeasures2018}.
A practical form of the predictability assumption is that there exists a constant $C>0$ such that $\observed(q)\leq C\predicted(q)$ for (almost every) $q\in\dspace$.
Conceptually, we interpret the predictability assumption as stating that we are able to predict all the observed data.
If the predictability assumption is satisfied, then as a consequence of $\updated$ being a density it follows that
\begin{linenomath*}
\begin{equation} 
	1 = \int_{\pspace} \initial(\param)\frac{\observed(Q(\param))}{\predicted(Q(\param))}d\initialP = \int_{\mathcal{D}} \frac{\observed(Q(\param))}{\predicted(Q(\param))}d\predictedP = 	\int_{\mathcal{D}} r(Q(\param))d\predictedP = \mathbb{E}(r) \label{eq:exp_r_deriv},
\end{equation}
\end{linenomath*}
where 
\begin{linenomath*}
\begin{equation} 
	r(Q(\param)):= \frac{\observed(Q(\param))}{\predicted(Q(\param))}.\label{eq:r}
\end{equation}
\end{linenomath*}
Thus, $\mathbb{E}(r)$ provides a convenient diagnostic for the computed MUD solution from a set of samples.
Specifically, if the sample average of $\mathbb{E}(r)$ deviates too far from unity, then this informs us of potential violations of the predictability assumption or other sources of error (such as density approximation error), in the formulation of the updated density. 
In such cases, further analysis is needed to determine the exact issue, but it nonetheless proves invaluable in determining the trustworthiness of the updated density and any statistical inferences drawn from this density.
We make extensive use of this diagnostic in Section~\ref{sec:PCA_examples}.

This diagnostic also helps to frame the special role of $\initial$ in the SIP as compared to the role of the prior density used in the Bayesian inverse problem that is discussed below.
Specifically, $\initial$ serves three roles in the definition and solution of the SIP that we emphasize in the list below.
\begin{itemize}
	\item $\initial$ represents an initial description of the {\em aleatoric} uncertainty on $\pspace$. 
	\item $\initial$ and the associated $\predicted$ allows us to construct a particular data-consistent solution in the form of $\updated$ that is {\em unique} up to the specification of $\initial$.
	\item $\updated$ differs from $\initial$ in $\pspace$ {\em only} in the directions for which $Q(\param)$ exhibits sensitivity. In other words, $\updated$ and $\initial$ have identical conditional probability structures on the manifolds defined by $Q^{-1}(q)$ for (almost every) $q\in\dspace$, which is evident by $r(Q(\lambda))$ being constant on such manifolds. We often refer to these manifolds as the {\em generalized contours} of the map $Q$. 
\end{itemize}

As discussed in the introduction, the major contributions and focus of this work are to utilize the data-consistent framework associated with the SIP to solve the PIP, which we define below.
\begin{definition}[Parameter Identification Problem (PIP)]\label{def:PIP}
Given a finite amount of (possibly noisy) data on a QoI map obtained for a fixed, but unknown, parameter $\paramref$, the PIP is to estimate $\paramref$. 
\end{definition}

The estimate of $\paramref$ we consider in this work is given by the maximal updated density (MUD) point defined as
\begin{linenomath*}
\begin{equation}\label{eq:mudpt_inital_defn}
	\mudpt := \argmax {\updated}(\param).
\end{equation}
\end{linenomath*}

\subsubsection{A Bayesian approach for solving the PIP}

We now develop a typical Bayesian approach for solving the PIP following the framework described in \cite{Stuart10}.
Recalling that $\paramref$ refers to the true parameter value, let $d$ denote the ``noisy'' data obtained on $Q(\paramref)$, which is often represented as
\begin{linenomath*}
\begin{equation*}
	d = Q(\paramref) + \xi,
\end{equation*}
\end{linenomath*}
where $\xi$ is a random variable used to model the measurement error that is often assumed to follow a Gaussian distribution.
Then, the data-likelihood function, often written as a conditional density, $\pi_\text{like}(d\, |\, \param)$, is formed.
This describes the differences in relative likelihoods that the data could have been generated from a particular $\param$.
Ideally, the largest values of $\pi_\text{like}(d\, | \, \param)$ occur whenever $\param$ is a ``good'' approximation to the true parameter $\paramref$.
The data-likelihood function is distinct from the observed density used in the data-consistent framework.

The next ingredient in the Bayesian framework is the specification of a prior density denoted by $\pi_\text{prior}(\param)$.
The prior describes the different relative likelihoods assumed for the true parameter before data are collected.
In other words, the prior represents an initial description of the {\em epistemic} uncertainty on $\pspace$.
This immediately distinguishes the role of the prior from the role of the {\em initial} density used in the data-consistent framework.

The posterior density (i.e., the solution to the Bayesian inverse problem) is given by a conditional density, denoted by $\pi_\text{post}(\param\, | \, d)$, proportional to the product of the prior and data-likelihood function.
In other words,
\begin{linenomath*}
\begin{equation*}
	\pi_\text{post}(\param\, | \, d) \propto \pi_\text{prior}(\param)\pi_\text{like}(d\, | \, \param)
\end{equation*}
\end{linenomath*}
This form of the density follows from Bayes' rule (not from the Disintegration Theorem as with the updated density).
The posterior can be interrogated to assess the difference in relative likelihoods of a fixed parameter given the observed data.
Subsequently, the posterior is often used to produce a ``best'' estimate of the true parameter.
For example, the maximum a posteriori (MAP) point is the parameter that maximizes the posterior density.

\begin{table}[htbp]
\centering
\begin{tabular}{|c|c|}
\hline
 & \\
$\displaystyle \updated(\param) = \initial(\param) \frac{\observed(Q(\param))}{\predicted(Q(\param))}
$
&
$
	\displaystyle \pi_{\text{post}}(\param\,|\,d) = \frac{\pi_{\text{prior}}(\param)\pi_\text{like}(d\,|\,\param)}{\int_{\pspace} \pi_\text{like}(d\, |\, \param)  \pi_{\text{prior}}(\param) d\mu_{\pspace}}
$
 \\ & \\ \hline
\end{tabular}
\caption{Updated density solving the data-consistent inverse problem (left) and posterior density solving the Bayesian inverse problem (right). The role of observable data, denoted by $d$, is made explicit in the posterior whereas the role is implicit in the updated in how it is used to construct the observed density, denoted by $\observed$.}
		\label{tab:dens_comparisons}
\end{table}

\subsection{Comparisons and an illustrative example}\label{sec:comparison_example}

Before we utilize an illustrative example to compare the posterior and updated densities (along with the associated MUD and MAP points), we find it useful to summarize these densities side-by-side in Table~\ref{tab:dens_comparisons} and comment on a few notable aspects not mentioned above.
Observe for the posterior density that the data-likelihood function appears in both the numerator and denominator.
In particular, the data-likelihood function informs the {normalizing constant}\footnote{The normalizing constant is commonly referred to as the evidence term.} in the denominator.
This is in contrast to the denominator of the updated density, which is given by the predicted density, which is in general not a constant, and can be constructed independently of $\observed$.

A practical implication of this difference is that the updated density only alters the structure of the initial density in what we refer to as the ``data-informed'' parameter directions.
Specifically, for a fixed $q\in\dspace$, let $C_q := \set{\param\in\pspace\, : \, Q(\param)=q}$, i.e., $C_q$ is a ``contour'' in parameter space.
Then, for any $\param\in C_q$, we immediately have $\updated(\param)=r(q)\initial(\param)$ where $r(q)$ is a fixed constant of proportionality for all $\param\in C_q$.
By contrast, while the posterior does not have to agree with the prior in any direction in parameter space, the prior does impact the structure of the posterior in all directions.

The previous paragraph is not\---and should not be interpreted as\---a criticism of the Bayesian inverse framework.
It is only meant to highlight that the data-consistent and Bayesian frameworks formulate and solve inverse UQ problems from different perspectives and with different (although at times seemingly compatible) assumptions.
Consequently, the solutions for an inverse problem formulated under either framework may differ significantly.
As the example (adopted from \cite{butlerCombiningPushForwardMeasures2018}) below demonstrates, this is true even if we arbitrarily force the inverse problems to appear as similar as possible.

Before diving into the comparison example, we quickly discuss the computational costs between computing the MUD and MAP point in this work. 
The initial/prior and observed/likelihood PDFs are either usually specified exactly as part of the setup for the SIP and PIP or have their forms determined from data separate from model evaluations.
However, for the density-based solutions to the SIP, we must generally approximate the predicted density.
We utilize random samples and non-parametric kernel density estimates (KDEs) to approximate the various pushforward densities as well as the updated and posterior densities since this is a straightforward method for implementation\footnote{See~\ref{app:A} for details on acquiring the datasets and code utilized for all examples in this manuscript.}.
Thus, the computational bottleneck in this work for both methods is the number of samples required to accurately estimate densities via a KDE, which is itself primarily dominated by the cost of evaluating these samples through our forward model (if it is expensive to evaluate).
As mentioned at the beginning of Section~\ref{sec:compare}, there are ``set-based'' solution methodologies (e.g., see \cite{butler2014measuretheoretic}), but these have issues with scaling to higher-dimensions.
Another approach suggested in \cite{marvin2018scalable} yields promising results for data-consistent inversion that alleviates the curse of dimensionality by adapting the Laplace approximation method that is commonly deployed in Bayesian settings.
Since the focus of this work is on the analysis of the MUD point and its utilization as an alternative to a MAP point, we leave to future work any analysis and comparison of alternative numerical approaches.

\subsubsection{Illustrative example: density comparison}

\begin{figure}[htbp]
\centering
   \includegraphics[width=0.49\linewidth]{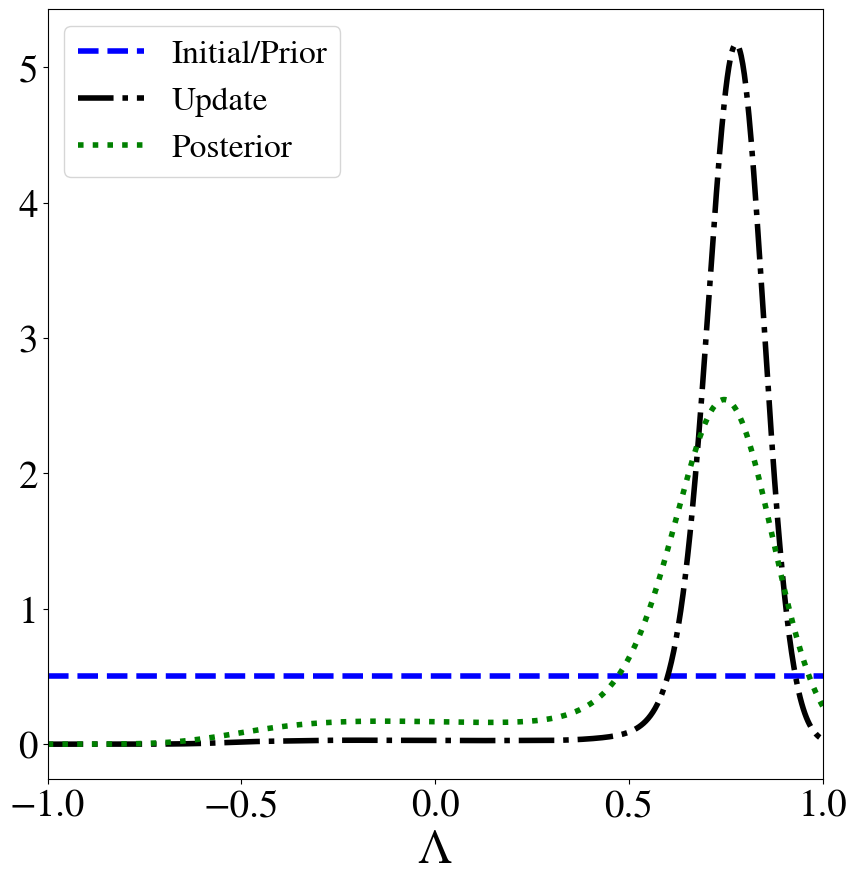}
   \includegraphics[width=0.49\linewidth]{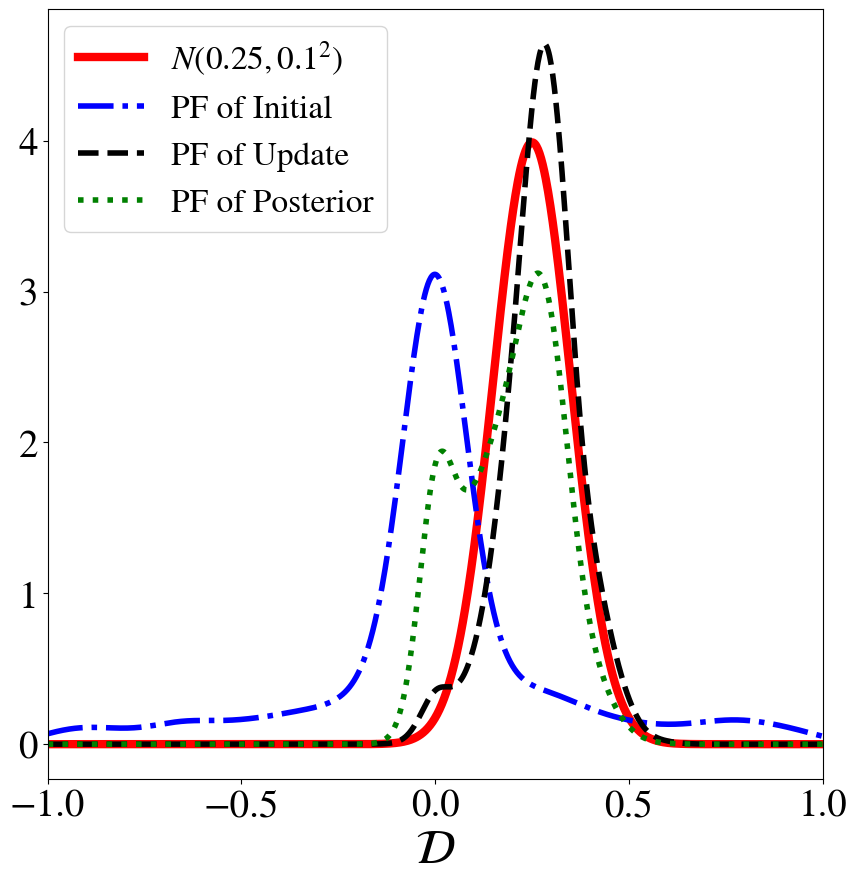}
 \caption{(Left) The initial/prior PDF $\initial$ (blue dashed curve), updated PDF $\updated$ (black dashed curve), and posterior PDF $\pi_\text{post}$ (green dashed-dotted curve) on $\pspace$.
 (Right) The push-forward (PF) of the initial/prior PDF $\predicted$ (blue dashed curve), observed/likelihood PDF (red solid curve), PF of the updated PDF $\updated$ (black dashed curve), and the PF of the posterior PDF $\pi_\text{post}$ (green dashed-dotted curve) for the QoI.}
 \label{fig:bayes-comparison}
\end{figure}

Suppose $\pspace = [-1,1]\subset\RR$ and $Q(\param)=\param^5$ so that $\dspace = [-1,1]$.
For the data-consistent framework, we assume $\initial\sim \mathcal{U}([-1,1])$ and $\observed\sim N(0.25,0.1^2)$.
The push-forwards of the initial, observed, and updated PDFs are shown in Figure~\ref{fig:bayes-comparison}.

For the Bayesian inverse problem, we assume $d\in \dspace$ with $d=Q(\paramref)+\xi$ where $\xi\sim N(0,0.1^2)$.
We then construct $\pi_{\text{post}}(\param \, |\, d)$ for this example assuming a uniform prior (to match the initial density) with an assumed observed value of $d=0.25$ so that the data-likelihood function matches the observed density.
The posterior and its push-forward are also shown in Figure~\ref{fig:bayes-comparison}.

While the updated and posterior densities in Figure~\ref{fig:bayes-comparison} share certain similarities (e.g., they are uni-modal with similar locations of the mode), they are otherwise visibly distinct.
The differences between these densities are made more evident by examining their push-forwards.
The push-forward of the updated density agrees well with the observed density, which is to be expected.
However, the push-forward of the posterior is bi-modal and does not match the observed density, which we recall is identical to the data-likelihood function in this case.

%

The takeaway is that each density is constructed to provide a solution to a {\em different} inverse problem.
The posterior density is intended to provide estimates of a true parameter value (i.e., the posterior is designed to solve a PIP) whereas the updated density is intended to quantitatively characterize natural variations in parameter values (i.e., the updated density solves a SIP).

\subsubsection{Illustrative example: MUD vs. MAP comparison}

To keep this comparison as clear as possible, we initially consider the exact MUD and MAP estimates associated with the exact updated and posterior densities. 
In other words, we avoid any approximation errors due to the use of estimated densities (for the MUD point) and techniques often utilized for approximating solutions to the optimization problems (for both the MUD and MAP points). 
Given the same setup as above, the exact MUD point obtained by utilizing the exact updated density is $\mudpt = 0.25^{1/5}$, which corresponds exactly to the value of $\paramref$ such that $Q(\paramref)=0.25$. 
Since the prior density is assumed uniform, the exact MAP point obtained by utilizing the exact posterior is identical to the MUD point.

It is worth noting that the MUD point in this example is invariant to the choice of initial density as long as the predictability assumption is not violated. 
In other words, the MUD point will always equal $\paramref$, i.e., it is always an unbiased estimate of $\paramref$ in this example.
However, if we altered the prior in this example to introduce a biased initial guess away from $\paramref$, the MAP point would no longer equal the MUD point because of the persistent bias present in the prior density. 

This is not to say that the MUD point is necessarily a better estimate of $\paramref$ compared to the MAP point.
One of the advantages of the MAP point is that the Bayesian formulation provides a natural quantification of uncertainty in the MAP point in terms of the variance in this estimate, which is a function of the amount of data utilized to produce it.
We discuss this in the next section.

\section{Data-Constructed QoI maps I: The Repeated Measurement Case}\label{sec:data-maps}
While the previous section establishes the MUD point as a potential alternative parameter estimate to the MAP point, we have yet to address the reduction of variance in MUD estimates as more data are included, which more naturally occurs for the MAP estimates.
The mechanism for reducing variance in MAP estimates is the data-likelihood function, which is often written as a product of conditional densities associated with the individual components of the data vector.
In other words, as more data are incorporated, the Bayesian approach is effectively increasing the dimension of the data space leading to a reduction of the variance in the subsequent MAP estimates.
This is not surprising since the Bayesian formulation is within the context of a PIP.
However, the MUD point is associated with the solution to a SIP for which a QoI map is defined and the dimension of the QoI data space is then fixed.
This suggests that we either incorporate data into the QoI map construction, or into the specification of the observed density in order to make the MUD point and updated covariance more useful for solving a PIP.
Below, we introduce some basic notation and return to the illustrative example from Section~\ref{sec:comparison_example} to motivate the utilization of data in the construction of a QoI map rather than in the specification of the observed density.

\subsection{Notation and an illustrative example}

Suppose there exists $m$ measurement devices for which repeated noisy data are obtained.
For each $1\leq j\leq m$, denote by $\mathcal{M}_j(\paramref)$ the $j$th measurement device, and denote by $N_j$ the number of noisy data obtained for $\mathcal{M}_j(\paramref)$.
Let $d_{j,i}$ denote the $i$th noisy datum obtained for the $j$th measurement where $1\leq i\leq N_j$.
We assume an unbiased additive error model for the measurement noise with independent identically distributed (i.i.d.) Gaussian errors so that
\begin{linenomath*}
\begin{equation}\label{eq:obs_data_error}
	d_{j,i} = \mathcal{M}_j(\paramref) + \xi_i, \ \xi_i\sim N(0,\sigma_j^2), \ \ 1\leq i\leq N_j.
\end{equation}
\end{linenomath*}
When a single measurement device is used to collect data, we often drop the index $j$.

\begin{figure}[htbp]
\centering
   \includegraphics[width=0.45\linewidth]{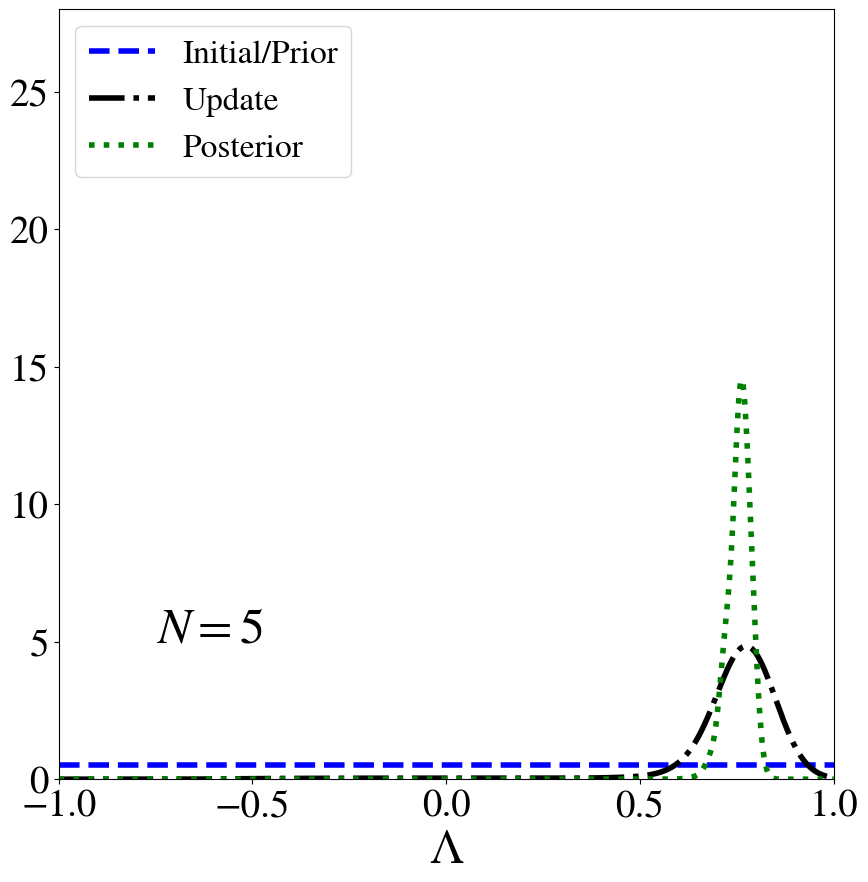}
   \includegraphics[width=0.45\linewidth]{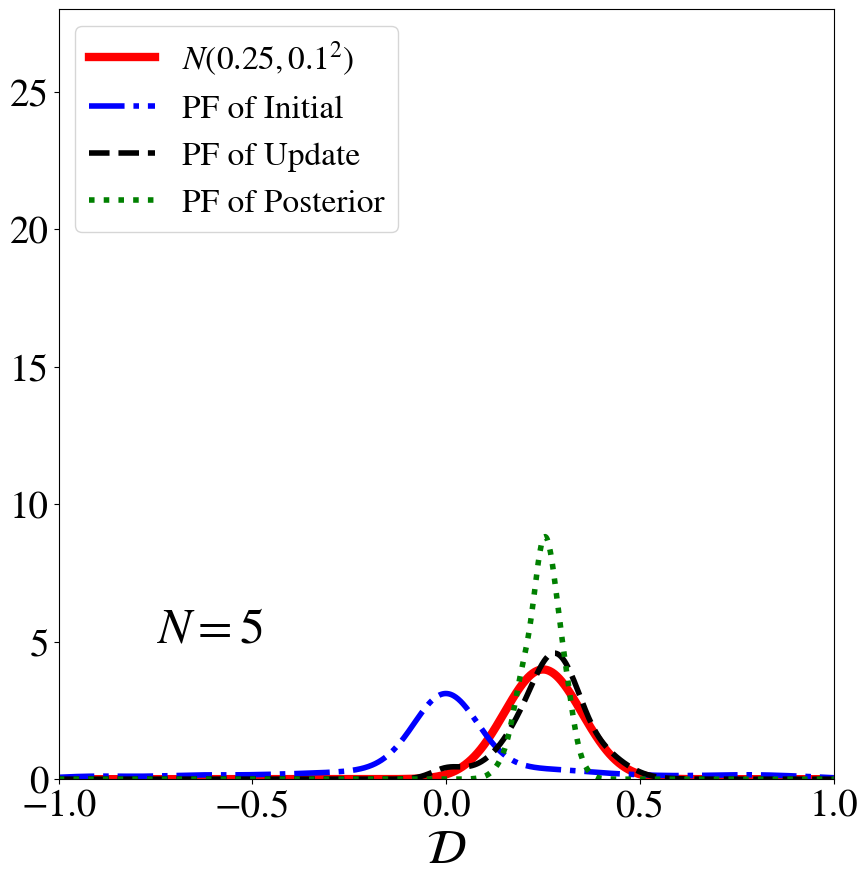}
   \includegraphics[width=0.45\linewidth]{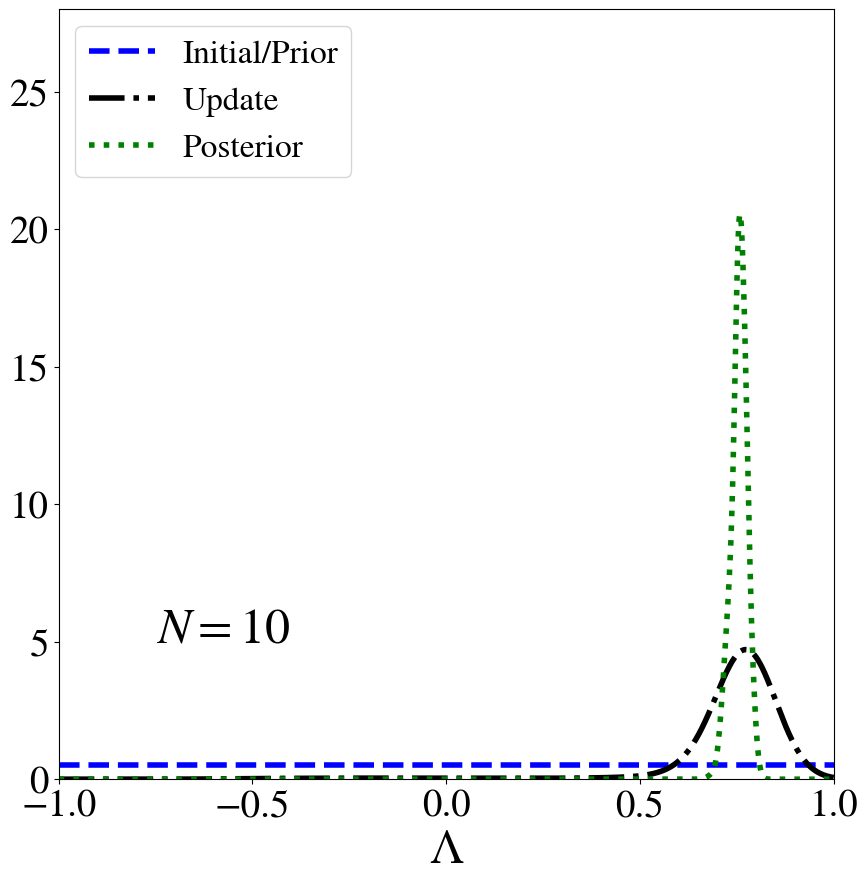}
   \includegraphics[width=0.45\linewidth]{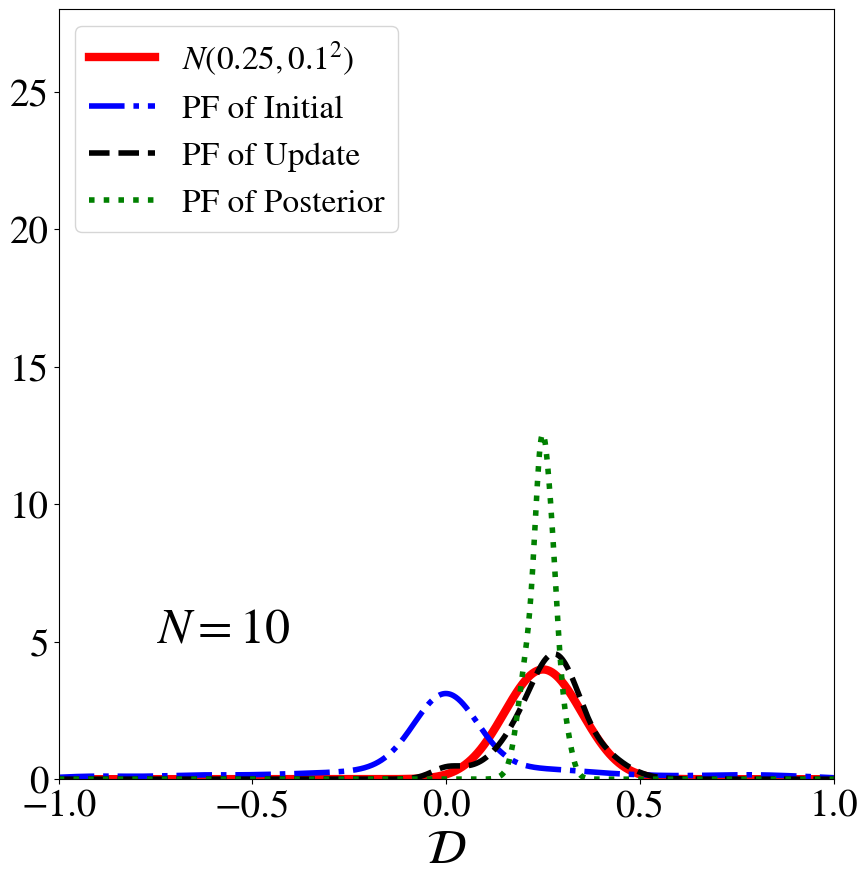}
   \includegraphics[width=0.45\linewidth]{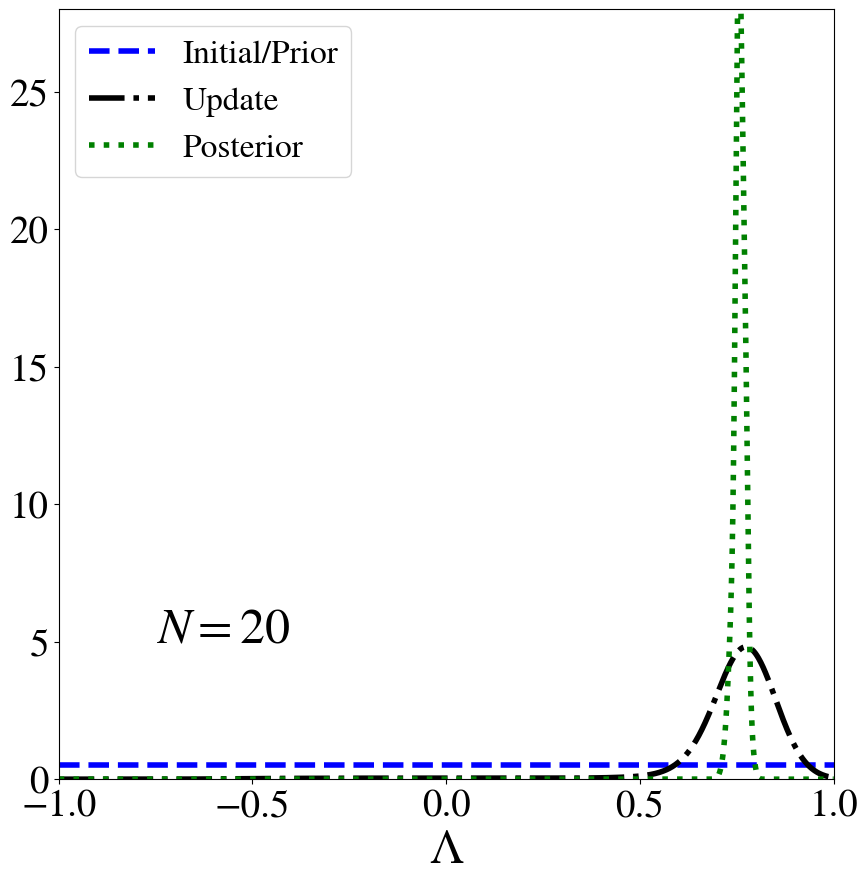}
   \includegraphics[width=0.45\linewidth]{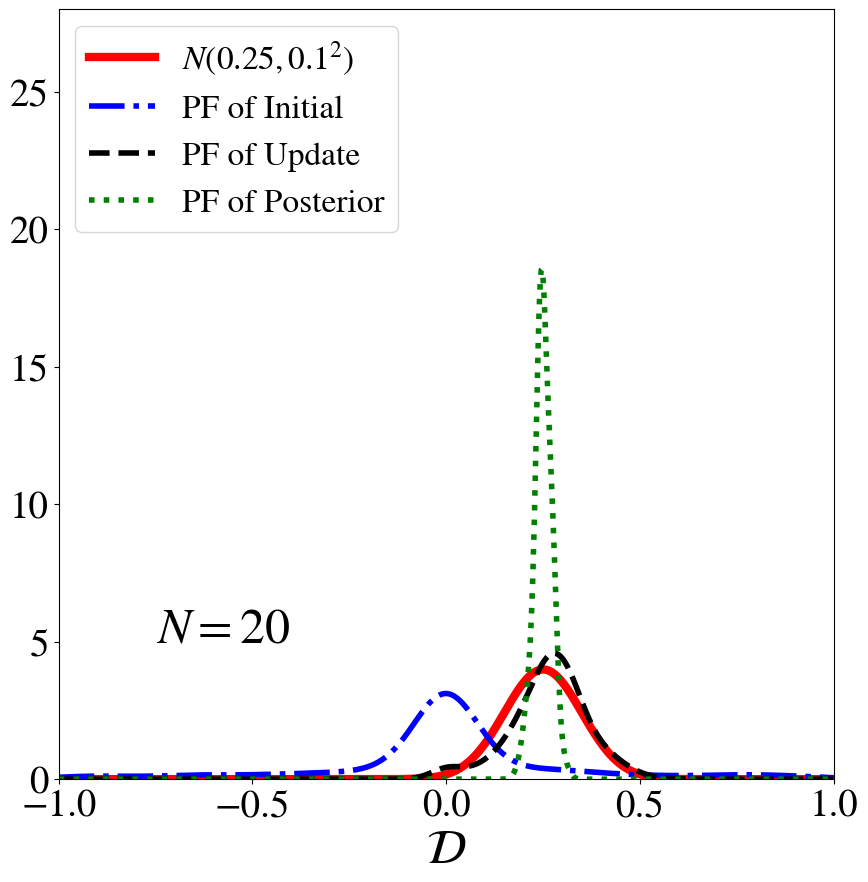}
 \caption{(Top to Bottom): $N=5, 10, \text{ and } 20$ samples are used to solve the inverse problems for comparison.
	  (Left) The initial/prior PDF $\initial$ (blue solid curve), updated PDF $\updated$ (black dashed curve), and posterior PDF $\pi_\text{post}$ (green dashed-dotted curve) on $\pspace$.
	  (Right) The push-forward (PF) of the initial/prior density (i.e., $\predicted$) (blue solid curve), observed/likelihood PDF (red solid curve), PF of the updated PDF $\updated$ (black dashed curve), and the PF of the posterior PDF $\pi_\text{post}$ (green dashed-dotted curve) for the QoI.}
 \label{fig:bayes-comparison-convergence}
\end{figure}

We now return to the example from Section~\ref{sec:comparison_example}.
Recall that for the Bayesian inverse problem, $Q(\paramref)=0.25$ and noisy data are given by $d=Q(\paramref)+\xi$ where $\xi\sim N(0,0.1^2)$.
The data are used to construct the data-likelihood function.
For the SIP, suppose that we use the sample mean of data to estimate the observed $N(0.25,0.1^2)$ distribution.
In other words, we use the sample mean of the data as an estimate for the mean of a normal distribution where the variance is assumed---as in the Bayesian inverse problem---to be $0.1^2$.
The observed density and data-likelihood become significantly different from one another as more data are collected.
The data-likelihood is in fact given by a product of normal densities evaluated at residuals of the data and the QoI map.
From this perspective, we interpret Figure~\ref{fig:bayes-comparison} as comparing the updated and posterior densities when we collect a single datum that happens to agree with the mean of the distribution. 
We now draw $N=5, 10, \text{ and } 20$ samples to form estimates of $\observed$ and the likelihood functions and show representative results in Figure~\ref{fig:bayes-comparison-convergence}.

For all values of $N$, the densities associated with the SIP all roughly agree.
In fact, the initial and predicted densities never change by construction. 
As more data are incorporated (i.e., as $N$ increases), the observed density estimates will stabilize as the variance in the mean of this density is reduced at a rate given by the central limit theorem.
Subsequently, the updated densities also stabilize and the MUD points are all generally in good agreement with $\paramref$.
However, we observe that the posterior becomes more peaked with a MAP point that converges to the MUD point as $N$ increases.
Thus, the uncertainty associated with the MAP point, i.e., the variance in the MAP point estimate for any given realization of $N$ data, is reduced as more data are incorporated.
This demonstrates that incorporating data into the formulation of the SIP to produce more stable estimates of the observed density will not necessarily reduce the variance in the MUD point estimates. 
We therefore seek an alternative approach to incorporating data into the SIP through the construction of the QoI map itself that will result in reduced variance in MUD estimates as more data are incorporated.

\subsection{The Weighted Mean Error Map}

We examine the data-likelihood function utilized in the Bayesian inverse problem and the concept of {\em statistical sufficiency} to motivate the data-constructed QoI map for MUD estimation.
For simplicity, we initially consider a single measurement device, drop the $j$ notation shown in~\eqref{eq:obs_data_error}, and let $\sigma$ denote the standard deviation of the measurement error. 
In the Bayesian inverse problem, the data-likelihood function formed by observing $N$ data points, $\set{d_i}_{i=1}^N$, has the form $ \pi_\text{like}(d\, |\, \param)$ where $d$ denotes an $N$-tuple with $i$th element given by $d_i$.
Assuming i.i.d. noise present in each datum, this joint data-likelihood function takes the form
\begin{linenomath*}
\begin{align} 
 \pi_\text{like}(d \mid \lambda) &= \prod_{i=1}^N \pi_\text{like}(d_i \mid \lambda) \nonumber \\
&\propto \prod_{i=1}^N \exp\left(-\frac{1}{2\sigma^2}(d_i-\mathcal{M}(\lambda)
)^2\right) \nonumber \\
&= \exp\left(-\frac{1}{2\sigma^2}\sum_{i=1}^n(d_i-\bar{d})^2\right)\exp\left(-\frac{1}{2(\sigma^2/N)}(\bar{d}-\mathcal{M}(\lambda))^2\right), \label{eq:fisher-neyman}
\end{align}
\end{linenomath*}
where $\bar{d}$ denotes the average of the $N$ elements of $d$. To go from the middle to last line above, we perform a classic re-factorization of the quadratic terms required by the Fisher-Neyman factorization theorem (see \cite{casella2002statistical} for details).
The Fisher-Neyman theorem states that $T(d):=\frac{1}{N}\sum d_i =\bar{d}$ is a sufficient statistic for $\lambda$ since the data-likelihood depends on the data only through the statistic $T$.
Thus, the likelihood principle suggests that any two sets of data with the same sample mean $\bar{d}$ will provide the same evidence or inference about $\lambda$.
If we wish to determine the variability in inferences about $\lambda$ due to differences in two separate collections of $N$ data points, it ``suffices'' to consider the variability in the statistic $\bar{d}$.

Motivated by this perspective, we first consider the data-constructed QoI map defined by the mean error (ME) and denoted by $Q_\text{ME}(\param)$, with $j$th component $Q_{\text{ME},j}(\param)$, given by 
\begin{linenomath*}
\begin{equation} 
	Q_{\text{ME},j}(\param) := \frac{1}{N_j} \sum_{i=1}^{N_j} \left(\mathcal{M}_j(\param)-d_{j,i}\right). \label{eq:qoi_ME}
\end{equation}
\end{linenomath*}
This choice of map leads to a predicted and observed distribution that both depend on the observed data characteristics, e.g., $Q_{\text{ME},j}(\paramref)\sim N(0,\sigma_j^2/N_j)$ which implies that $\observed$ is defined by a $N(\mathbf{0}_{m\times 1},\observedCov)$ distribution where $\observedCov$ is a diagonal matrix with $j$th diagonal element given by $\sigma_j^2/N_j$. 
As shown below, by applying some scaling transformations to the components of $Q_{\text{ME}}$, we can remove the dependence of each component of the observed distribution on the characteristics of the observed data, which simplifies the practical implementation.

The weighted mean error (WME) map, denoted by $Q_\text{WME}(\param)$ with $j$th component $Q_{\text{WME},j}(\param)$, is given by
\begin{linenomath*}
\begin{equation} 
	Q_{\text{WME},j}(\param) := \frac{1}{\sqrt{N_j}} \sum_{i=1}^{N_j} \frac{\mathcal{M}_j(\param)-d_{j,i}}{\sigma_j}. \label{eq:qoi_WME}
\end{equation}
\end{linenomath*}
The rationale for creating a WME map as opposed to the ME map is found by substitution of \eqref{eq:obs_data_error} into \eqref{eq:qoi_WME}.
Specifically, by this substitution followed by rationalizing the denominator, $Q_{\text{WME},j}(\paramref)$ is identified as the sample average of $N_j$ random draws from an i.i.d.~$N(0,N_j)$ distribution.
By assumption, all of the observed data are generated according to the true parameter vector given by $\paramref$ in \eqref{eq:obs_data_error}.
Subsequently, each component of $Q_\text{WME}(\paramref)$ is a random draw from an $N(0,1)$ distribution.
Therefore, with this choice of data-defined QoI map, we specify $\observed$ as a $N(\mathbf{0}_{m\times 1},\mathbf{I}_{m\times m})$ distribution.
In other words, this method of incorporating data into the SIP produces an observed distribution which is stationary with respect to the number of measurements used.
While this simplifies the algorithmic implementation of this approach (the observed distribution is fixed with respect to the data), as we show in Section~\ref{subsec:unifying-perspective}, the variance in MUD estimates is a decreasing function of the number of observed data.

\subsection{WME and the Predictability Assumption}\label{sec:MUD_analysis}

Since $\mudpt$ requires maximizing $\updated$, in this section we explore the predictability assumption for the WME map under certain assumptions which is sufficient for guaranteeing the existence and uniqueness of $\updated$ for a given $\initial$ and a QoI map.
This also serves to motivate Section~\ref{sec:estimation} exploring the linear Gaussian theory of existence and uniqueness of MUD points along with convergence of MUD estimates obtained by linear WME maps. 
In Section~\ref{sec:PCA}, we present an alternative to the WME map based on a principal component analysis (PCA) for situations where repeated data are not obtained on a fixed set of measurements. 

We begin by assuming $\initial \propto \mathcal{N}(\initialMean, \initialCov)$ with non-degenerate $\initialCov$.
For each $1\leq j\leq m$, assume that the observable maps  are linear maps of $\param$ (affine maps require minor modifications to the results below).
We denote these by $M_j$ instead of $\mathcal{M}_j$ to emphasize the linear assumption of these measurement maps.
For notational convenience below, assume that $M_j$ is written explicitly as a $1\times p$ row vector and that the $m$-vectors form a linearly independent set.
Then, it is possible to rewrite $Q_{\text{WME}}(\param)$ as
\begin{linenomath*}
\begin{equation} 
	Q_\text{WME}(\param) = A(\mathbf{N})\param + b(\mathbf{N}), \label{eq:q_wme_linear}
\end{equation}
\end{linenomath*}
where the $j$th component of $\mathbf{N}\in\RR^m$ is given by $N_j$ and the $j$th row of $A(\mathbf{N})\in\RR^{m\times p}$ is given by
\begin{linenomath*}
\begin{equation}
	\frac{1}{\sqrt{N_j}} \sum_{i=1}^{N_j} \frac{M_j}{\sigma_j} = \frac{\sqrt{N_j}}{\sigma_j}M_j,
\end{equation}
\end{linenomath*}
and the bias vector, $\mathbf{b}(\mathbf{N})\in\RR^m$, is defined by the data, with $j$th component, denoted by $\mathbf{b}_j$, given by
\begin{linenomath*}
\begin{equation}
	\mathbf{b}_j(\mathbf{N}) = -\frac{1}{\sqrt{N_j}} \sum_{i=1}^{N_j} \frac{d_{j,i}}{\sigma_j}.
\end{equation}
\end{linenomath*}
Since $A(\mathbf{N})\initialCov A(\mathbf{N})^\top$ defines a predicted covariance, and the observed covariance is the identity map, the predictability assumption is immediately satisfied if each diagonal component of the predicted covariance is significantly greater than $1$.

The off-diagonal components of the predicted covariance (which are dictated in large part by the structure of the measurement operators $M_j$) dictate how much larger than $1$ each diagonal component must be to ensure the predictability assumption holds.
However, we demonstrate below that this will happen once a minimum number of data points $N_\text{min}$ are obtained for each measurement.

First, observe that the $j$th diagonal component of the predicted covariance matrix is given by the predicted variance associated with using the scalar-valued map $Q_{\text{WME},j}$.
Then, the associated predicted variance is given by
\begin{linenomath*}
	\begin{equation}\label{eq:wme_pred_cov}
	\frac{N_j}{\sigma_j^2} M_j\initialCov M_j^\top
\end{equation}
\end{linenomath*}
Since $\initialCov$ is assumed to be non-degenerative and $M_j$ is a non-trivial row vector, this predicted variance grows linearly with $N_j$.
In other words, the $j$th diagonal component of the predicted covariance has the form $\beta_j N_j$ for some $\beta_j>0$.
Therefore, for each $1 \leq j \leq N$, there exists $N_{min,j}$ such that $N_j \geq N_{min,j}$ guarantees that the $j$th diagonal component is sufficiently large so that the smallest eigenvalue of the predicted covariance is larger than 1.
Finally, we observe that the predicted covariance is inversely proportional to the measurement noise present in~\eqref{eq:wme_pred_cov}, which indicates that more data points $N_j$ are required for measurement devices with greater noise.

\section{Linear Theory of MUD points}\label{sec:estimation}

In this section, we assume linear (or affine) QoI maps with Gaussian distributions, which are often used in the UQ literature to provide a common framework for comparing methods and their solutions.

This section is structured into several subsections to help focus the interpretations and results.
In Section~\ref{subsec:Motivation}, we present some useful details, notation, and terminology used for this comparison framework.
To build intuition, we compare both the MUD and MAP points in Section~\ref{subsec:low-d-example} using a low-dimensional example.
A unifying perspective is provided for affine maps in Section~\ref{subsec:unifying-perspective} along with derivations of closed form expressions for the MUD and MAP points in this comparison framework.
These results are summarized in a theorem of existence and uniqueness of the MUD point in this comparison framework.
This is followed by a corollary involving the convergence of MUD points obtained by WME maps and a brief analysis in Section~\ref{ex:spectral-decay} of the spectral properties of the updated covariance for a WME map as more data are incorporated into the map.
We then provide some higher-dimensional performance comparisons of MUD, MAP, and least squares estimates in Section~\ref{sec:high-dim-linear-example}. 

\subsection{Problem formulation and assumptions}\label{subsec:Motivation}

Let $\norm{\mathbf{x}}_C^2 := (\mathbf{x}, \mathbf{x})_C = \mathbf{x}^T C \mathbf{x} $ denote the square of the induced norm associated with a positive-definite operator $C : \RR^k \to \RR^k$ and the usual (Euclidean) inner product.
In what follows, the inverse covariances associated with non-degenerative multivariate Gaussian distributions will play the role of $C$.

Suppose that the initial and prior densities are both given by the same $\mathcal{N}(\param_0, \initialCov)$ distribution.
Additionally, suppose the map $Q$ is linear and that the data-likelihood and observed densities are both given by the same $\mathcal{N}(\observedMean, \observedCov)$ distribution.

The linearity of $Q$ implies that $Q(\param)=A\param$ for some $A\in\RR^{m\times p}$, and that the predicted density follows a $\mathcal{N}(Q(\param_0), \predictedCov)$ distribution where
\begin{linenomath*}
\begin{equation}\label{eq:predictCov}
	\predictedCov := A\initialCov A^\top.
\end{equation}
\end{linenomath*}
While it is not technically necessary to ensure that the predictability assumption holds (i.e, that $\updated$ is in fact a density) in order to formally define a MUD point using~\eqref{eq:mudpt_inital_defn}, it is useful when discussing certain theoretical results as shown in Section~\ref{subsec:unifying-perspective}.
Conceptually, the predictability assumption holds when the predicted variance is larger in all directions than the observed variance.
Mathematically, this occurs when the smallest eigenvalue value of $\predictedCov$ is larger than the largest eigenvalue value of $\observedCov$.
As discussed in Section~\ref{sec:data-maps}, this condition holds once a sufficient amount of data are observed for a WME map.

When $\observedCov$ is non-degenerative (i.e., the smallest eigenvalue is positive), the predictability assumption can always be satisfied if $m\leq p$ and $A$ is full rank by choosing $\initialCov$ to have sufficiently large eigenvalues (i.e., if we choose initial variances to be sufficiently large).
For clarity in the theoretical presentation of this section, we assume these conditions are met so that the predictability assumption holds and $\updated$ defines a density.
However, in Section~\ref{sec:high-dim-linear-example}, we still compute the formal MUD point for a high-dimensional example involving rank-deficient $A$ to demonstrate the overall usefulness of the MUD point even in situations where $\updated$ may fail to be an actual density.

\begin{table}
\centering
\begin{tabular}{|c|c|}
\hline
  & \\
  Tikhonov & $T(\param):=\norm{Q(\param)-\observedMean}_{\observedCov^{-1}}^2 +
      \norm{\param-\initialMean}_{\initialCov^{-1}}^2$
  \\ & \\ \hline & \\
  Data-Consistent & $J(\param):=T(\param) - \norm{Q(\param)-Q(\initialMean)}_{\predictedCov^{-1}}^2$
  \\ & \\
  \hline
\end{tabular}
\caption{Functionals to minimize to obtain $\param$ that maximizes  the updated PDF (bottom) and the Bayesian posterior PDF (top).
Here, $T(\param)$ is the typical functional often associated with Tikhonov regularization, and the $J(\param)$ has an additional term subtracted from $T(\param)$ coming from the predicted density that serves as ``un-regularization'' in data--informed directions.}
  \label{tab:func_comparisons}
\end{table}

With these assumptions, the parameters that maximize the posterior and updated densities are described as the arguments that minimize certain quadratic functionals.
Table~\ref{tab:func_comparisons} presents a scaling of these functionals defined by the negative logarithm of the associated posterior and updated densities.
Note that the functional, $T(\param)$, obtained from the posterior density is immediately identified as the typical functional used in Tikhonov regularization \citep{tarantola2005inverse}. The data-mismatch term given by
\begin{linenomath*}
\begin{equation*}
	\norm{Q(\param)-\observedMean}_{\observedCov^{-1}}^2
\end{equation*}
\end{linenomath*}
comes from the data-likelihood/observed density whereas the regularization term defined by
\begin{linenomath*}
\begin{equation*}
	 \norm{\param-\initialMean}_{\initialCov^{-1}}^2
\end{equation*}
\end{linenomath*}
comes from the prior/initial density.
We refer to this term as the Tikhonov regularization term.

By comparison, the functional, $J(\param)$, obtained from the updated density, is written as a modification of $T(\param)$ where the subtraction of
\begin{linenomath*}
\begin{equation*}
	 \norm{Q(\param)-Q(\initialMean)}_{\predictedCov^{-1}}^2
\end{equation*}
\end{linenomath*}
comes from the predicted density.

\subsection{A Low-Dimensional Example}\label{subsec:low-d-example}
To build intuition around the fundamental differences of MAP and MUD points beyond what is discussed in Section~\ref{sec:compare}, we consider an example where the linear QoI map is defined by $A=\mat{cc}{1 & 1}$, i.e., the parameter space is 2-dimensional while the data space is 1-dimensional.

In this example, the parameters in the initial and observed densities are given by
\begin{linenomath*}
\begin{equation*}
	\initialMean = \mat{cc}{0.25 & 0.25}^\top, \initialCov = \mat{cc}{1 & -0.25 \\ -0.25 & 0.5}, \ \observedMean=1, \text{ and } \observedCov = \mat{c}{0.25}.
\end{equation*}
\end{linenomath*}

\begin{figure}[htb]
  \centering
  \begin{tabular}{|ccc|}
    \hline
      \subf{\includegraphics[width=0.3\linewidth]{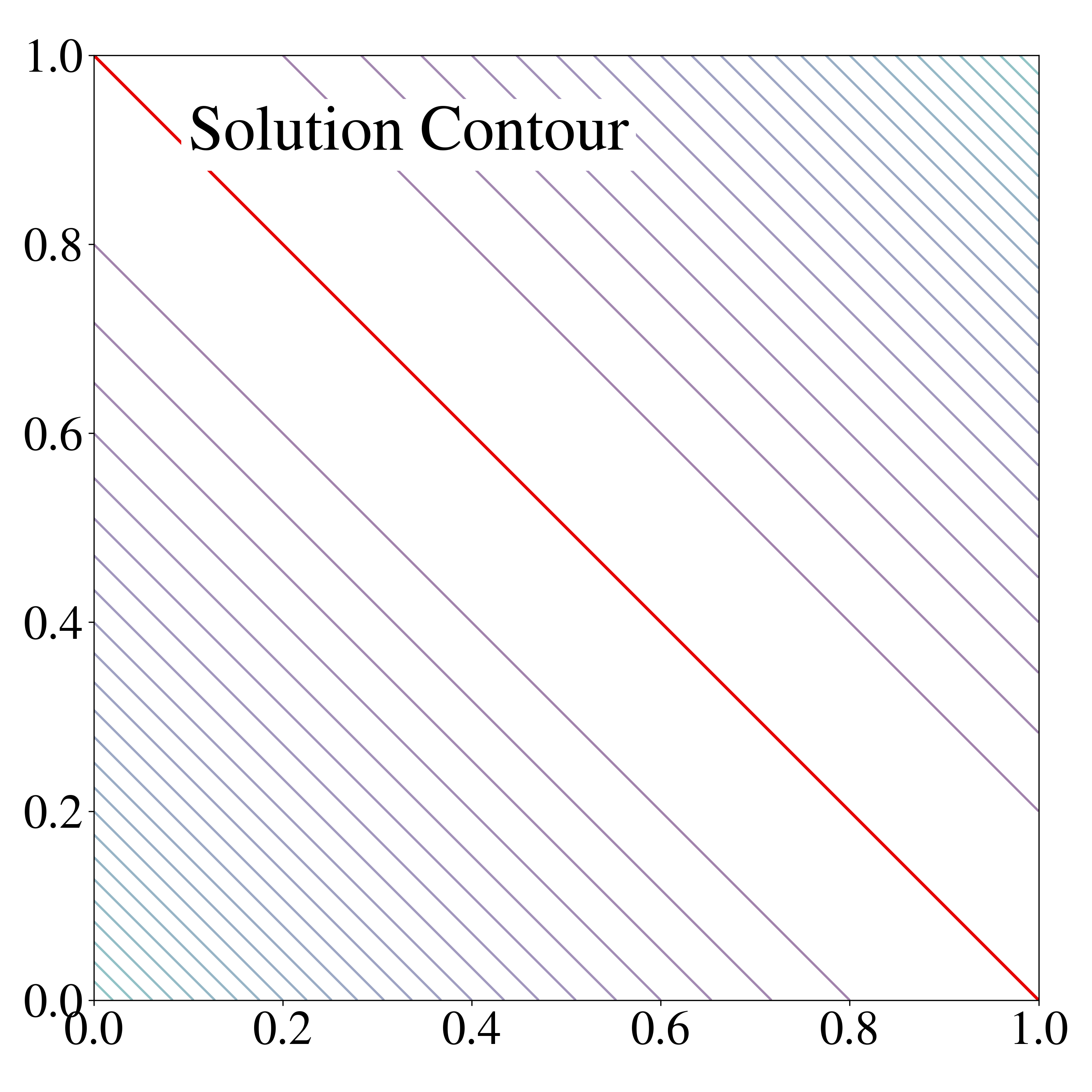}}
      {Data Mismatch}
    &
      \subf{\includegraphics[width=0.3\linewidth]{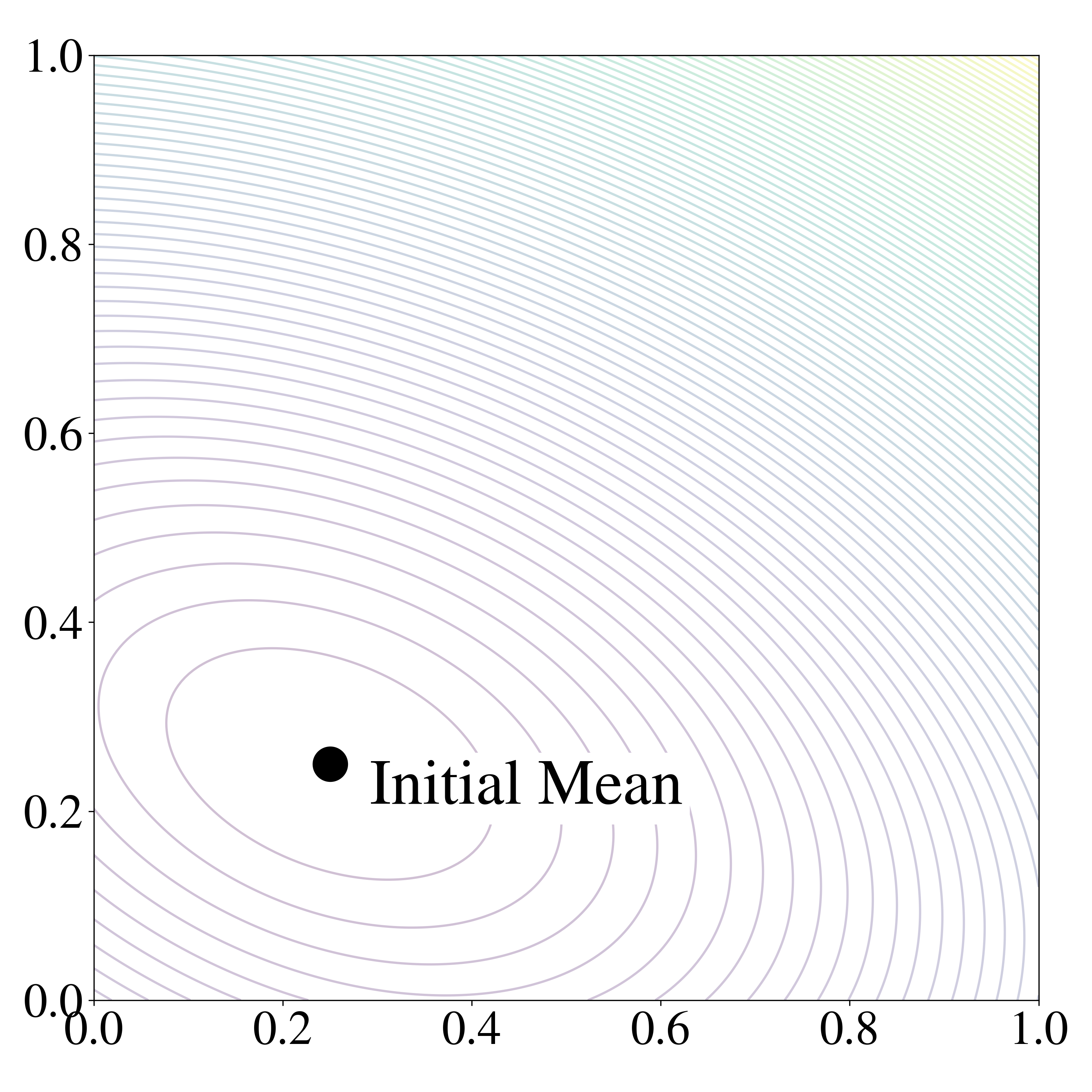}}
      {Regularization}
    &
      \subf{\includegraphics[width=0.3\linewidth]{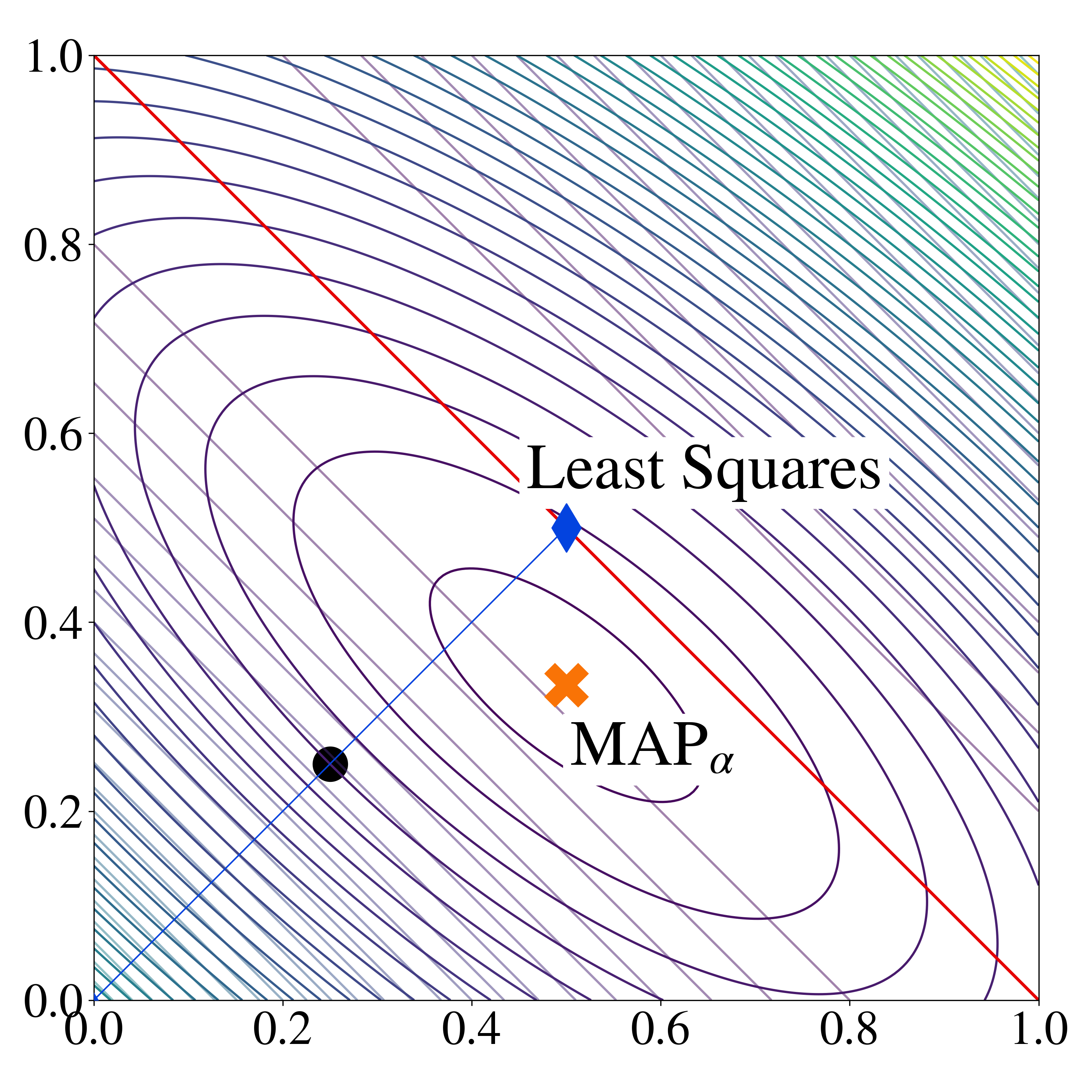}}
      {Bayesian Posterior}
    \\
    \hline
      \subf{\includegraphics[width=0.3\linewidth]{figures/contours/data_mismatch_contour.png}}
      {Data Mismatch}
    &
      \subf{\includegraphics[width=0.3\linewidth]{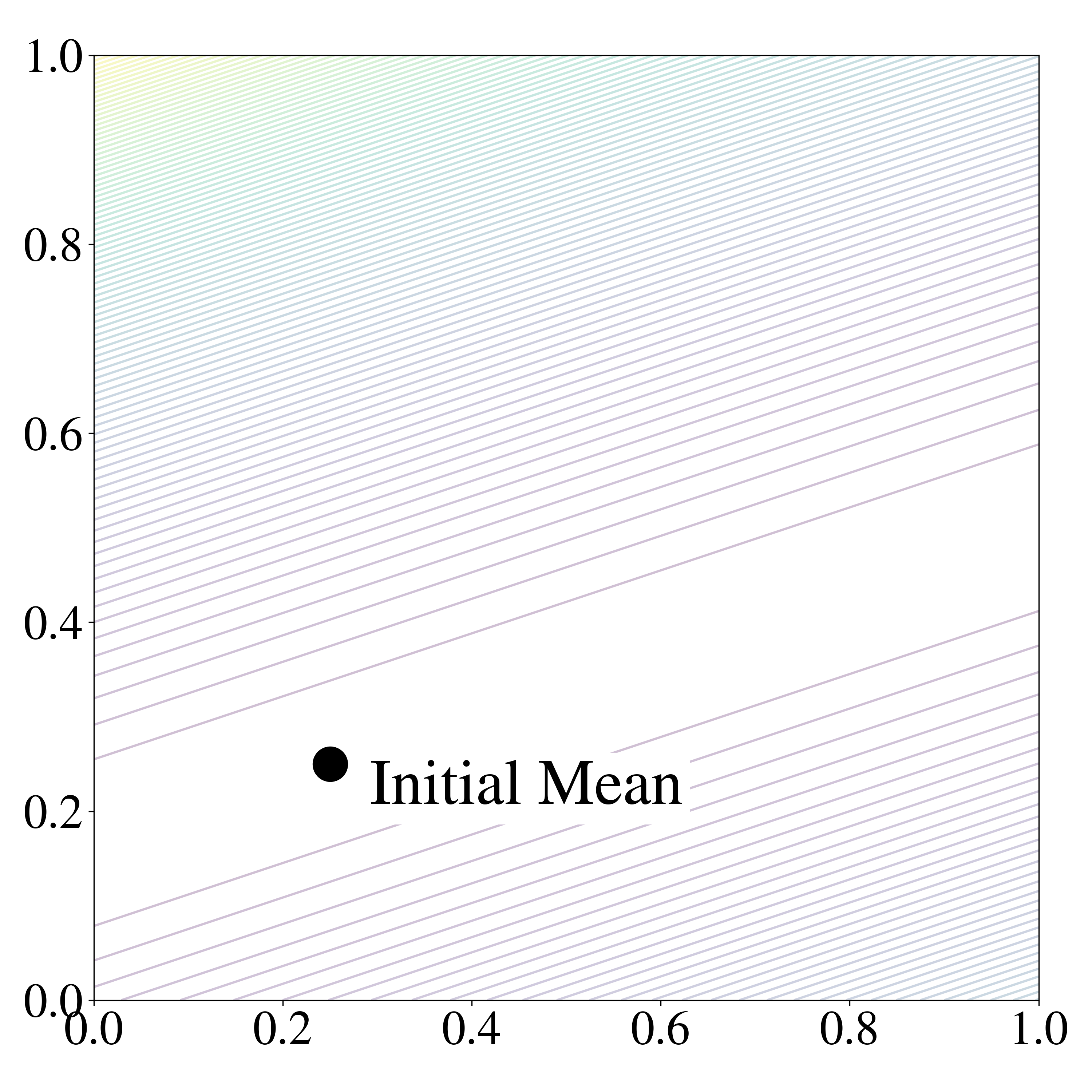}}
      {Modified Regularization}
    &
      \subf{\includegraphics[width=0.3\linewidth]{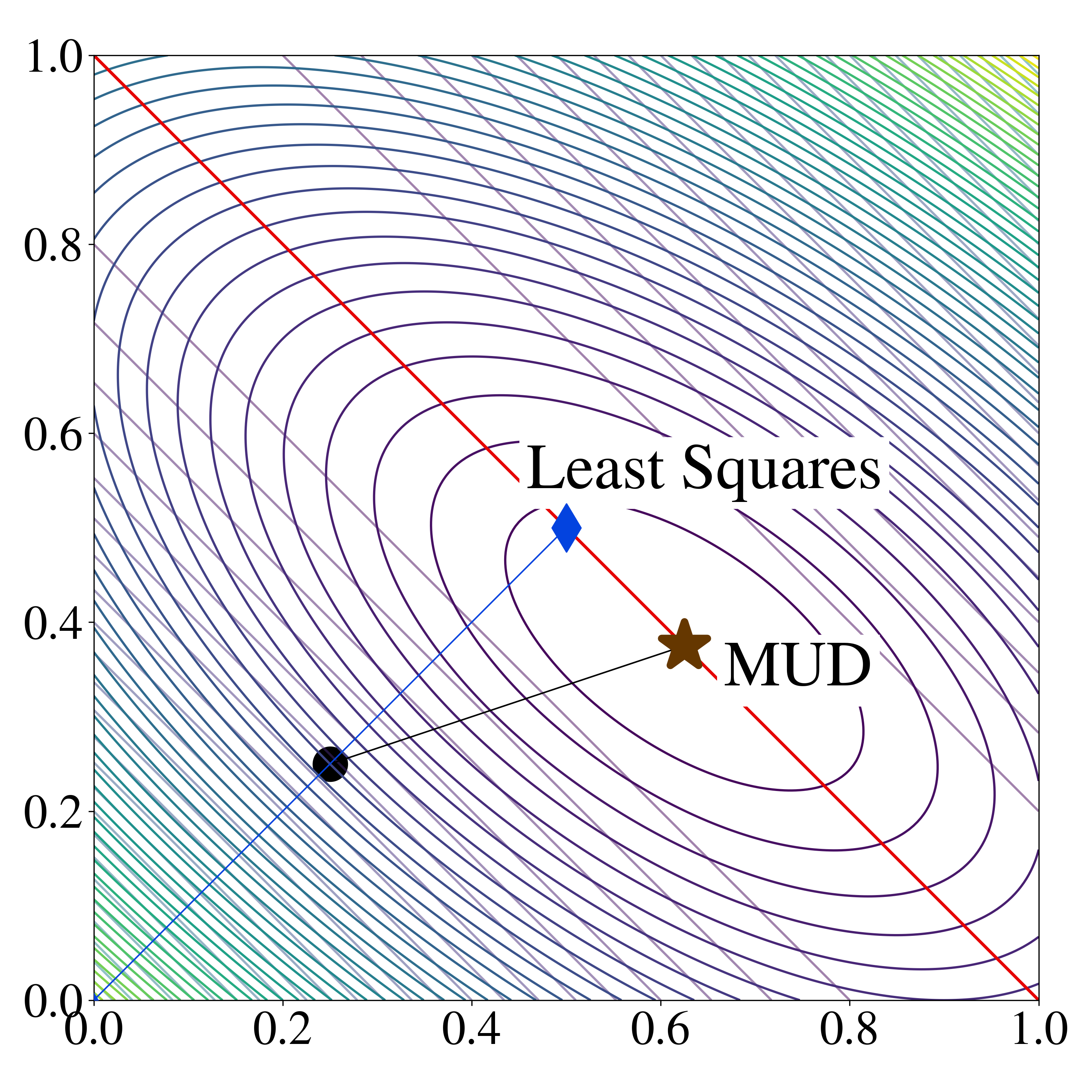}}
       {Updated Density}
    \\
    \hline
  \end{tabular}
  \caption{Gaussian data mismatch over 2-D parameter space for a 2-to-1 linear map (left plots). Gaussian initial/prior lead to different regularization terms associated with updated/Bayesian PDFs (middle plots), which lead to different optimization functions (right plots) and parameter estimates that produce maximum PDF values for update/Bayesian PDF (red dot in right plots).}
  \label{fig:regularization}
\end{figure}

The top row of Figure~\ref{fig:regularization} shows contour plots in the parameter space for the data-mismatch term (left), Tikhonov regularization term (middle), as well as the functional $T(\param)$ (right).
Conceptually, the regularization term is a radially symmetric function that penalizes parameters that are far away from the initial mean.

The bottom row of Figure~\ref{fig:regularization} shows contour plots in the parameter space for the data-mismatch term (left), modified regularization term (middle), as well as the functional $J(\param)$ (right).
Here, we see that the modified regularization term only penalizes the movement of parameters in certain directions away from the initial parameter mean.


\subsection{Existence and Uniqueness of MUD points}\label{subsec:unifying-perspective}

Assume that the QoI map, $Q$, now takes the slightly more general form
\begin{linenomath*}
\begin{equation}\label{eq:qoi_map}
Q(\param) = A \param + \mathbf{b}
\end{equation}
\end{linenomath*}
where $\mathbf{b} \in \RR^m$ may be viewed as a bias in the QoI map.
The inclusion of this term makes it possible to draw conclusions involving the data-constructed QoI maps presented in Section~\ref{sec:data-maps}.
Using the same Gaussian distribution assumptions as described in Section~\ref{subsec:Motivation}, we again identify the MAP and MUD points as the values that minimize the functionals $T(\param)$ and $J(\param)$, respectively, shown in Table~\ref{tab:func_comparisons}.

The posterior covariance is formally given by
\begin{linenomath*}
\begin{equation}\label{eq:map_cov}
\Sigma_\text{post} := ( A^\top \observedCov^{-1} A + \initialCov^{-1} )^{-1}.
\end{equation}
\end{linenomath*}
Applying the Woodbury matrix identity and~\eqref{eq:predictCov}, we rewrite the posterior covariance as
\begin{linenomath*}
\begin{equation}\label{eq:map_cov_analytical}
\Sigma_\text{post} = \initialCov - \initialCov A^\top \left[\predictedCov + \observedCov\right]^{-1} A \initialCov,
\end{equation}
\end{linenomath*}
which allows us to interpret $\Sigma_\text{post}$ as a rank $m$ correction (or update) of $\initialCov$.
Note that $\predictedCov + \observedCov$ is invertible because it is the sum of two symmetric positive definite matrices.
With either version of $\Sigma_\text{post}$ given above, we rewrite the closed form expression for the MAP point given in \cite{tarantola2005inverse} as
\begin{linenomath*}
\begin{equation}\label{eq:map-point-analytical}
\param^{\text{MAP}} = \param_0 + \Sigma_\text{post} A^\top \observedCov^{-1} (\observedMean - b - A\param_0).
\end{equation}
\end{linenomath*}

We can arrive at a similar expression for the $\lambda^{MUD}$ point by first deriving the updated covariance (see ~\ref{app:mud_deriv} for details), which is of the form
\begin{linenomath*}
\begin{equation}\label{eq:updatedCov_body}
	\updatedCov = \initialCov - \initialCov A^\top \predictedCov^{-1}\left[\predictedCov-\observedCov\right]\predictedCov^{-1}A\initialCov,
\end{equation}
\end{linenomath*}
which leads to
\begin{linenomath*}
\begin{equation}\label{eq:mud-point-analytical-final}
	\mudpt = \param_0 + \initialCov A^\top \predictedCov^{-1}(\observedMean - b - A\param_0).
\end{equation}
\end{linenomath*}

Comparing~\eqref{eq:mud-point-analytical-final} to~\eqref{eq:map-point-analytical}, we
see that the MUD point does not depend on the observed covariance whereas the MAP point does.
Moreover, applying $Q$ to \eqref{eq:mud-point-analytical-final} and substituting accordingly reveals that $Q(\mudpt) = \observedMean$.

Overall, this motivates the MUD point as an \emph{alternative parameter estimate} with predictive accuracy and properties directly correlated to the relationship between $\observedMean$ and the true signal for which noisy data are generated.
We summarize the above results in the following theorem stating the existence and uniqueness of a MUD point for the linear Gaussian case.

\begin{theorem}\label{thm:MUD_existence_uniqueness}
Suppose  $Q(\param)=A\param+b$ for some full rank $A\in\RR^{m\times p}$ with $m\leq p$ and $b\in\RR^m$.
If $\initial \sim N(\param_0,\initialCov)$, $\observed\sim N(\observedMean,\observedCov)$, and the predictability assumption holds, then
\begin{enumerate}[(a)]
\item There exists a unique parameter, denoted by $\mudpt$, that maximizes $\updated$.
\item $Q(\mudpt) = \observedMean$.
\item If $d=p$, $\mudpt$ is given by $A^{-1}$. If $m<p$, $\mudpt$ is given by~\eqref{eq:mud-point-analytical-final} and the covariance associated with this point is given by~\eqref{eq:updatedCov_final}.
\end{enumerate}
\end{theorem}

Recalling the discussion of Section~\ref{sec:MUD_analysis}, the following result is an immediate consequence of Theorem~\ref{thm:MUD_existence_uniqueness} and equation~\eqref{eq:mud-point-analytical-final},

\begin{corollary}\label{cor:MUD_wme}
If $\initial \sim N(\param_0,\initialCov)$ and data are obtained for $m$ linearly independent measurements on $\pspace$ with an additive noise model with i.i.d. Gaussian noise for each measurement, then
\begin{enumerate}[(a)]
\item There exists a minimum number of data points obtained for each of the measurements such that there exists a unique $\mudpt$ and $Q_\text{WME}(\mudpt) = 0$.
\item The variance in the $\mudpt$ estimate in $m$ directions of the parameter space decreases at a rate proportional to the number of data points used for each of the measurements, and is inversely proportional to the magnitude of the Gaussian noise, a relationship expressed in equation ~\eqref{eq:wme_pred_cov}.
\end{enumerate}
\end{corollary}

\subsection{Spectral Properties of the Updated Covariance for a WME map}\label{ex:spectral-decay}
\begin{figure}[htbp]
\centering
   \includegraphics[width=0.6\linewidth]{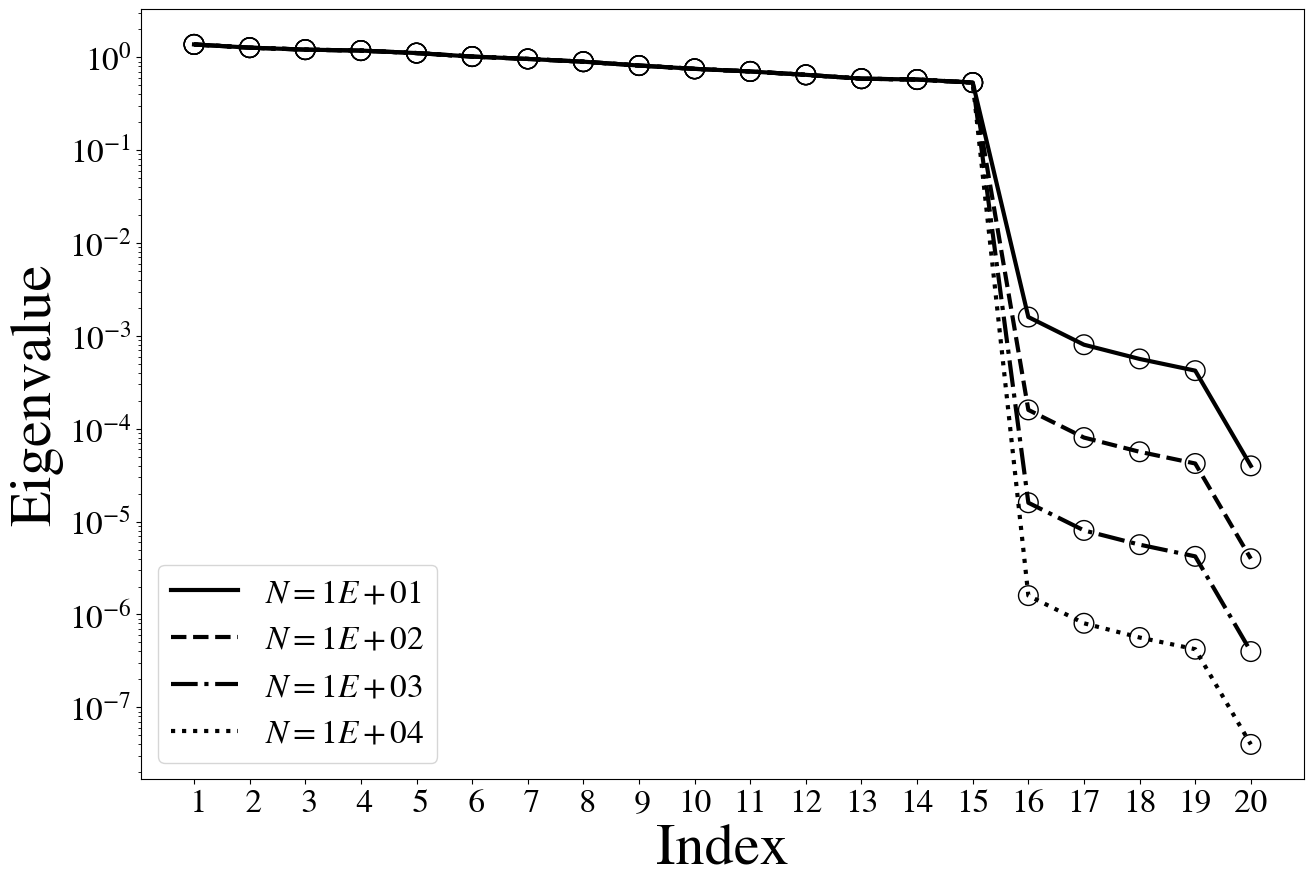}
 \caption{
	 Ranked eigenvalues of the updated covariance are shown for ${N=10, 100, 1000,}$ and $10000$ measurements, plotted by index as a line for each $N$.
   Fifteen of them show no sensitivity to $N$, and all are of $\mathcal{O}(1)$.
   With more measurements, the five eigenvalues (corresponding to the orthogonal null-space of the QoI operator) go to zero, losing approximately an order of magnitude for each increase in $N$.
 }
   \label{fig:linear-eigenvalues}
\end{figure}

We illustrate the result in Corollary~\ref{cor:MUD_wme} with an example that highlights how the number of distinct measurements used to form $Q$ impacts the spectrum of the updated covariance.
Consider a randomly generated linear operator $M$ of dimension $5 \times 20$ with components sampled from a $N(0,1)$ distribution.
This $M$ defines five randomly constructed, but geometrically distinct (i.e., linearly independent), measurements.
We then construct the QoI by simulating repeated measurements $N=10, 100, 1000, 10000$ (for each measurement) polluted by additive Gaussian noise with $\sigma=0.1$.
We then compute the updated covariance using the analytical expression \eqref{eq:updatedCov_final} and perform a singular value decomposition to obtain the $20$ sorted eigenvalues, which we plot for each $N$ in Figure~\ref{fig:linear-eigenvalues}.

Observe that five of the twenty eigenvalues are several orders of magnitude smaller than the rest, which corresponds to the output dimension of $Q$.
These correspond to the five directions informed by $Q$ given by the associated eigenvectors.
Furthermore, as observed in Fig~\ref{fig:linear-eigenvalues}, the gap between the uninformative and informative directions is directly proportional to $N$; for each ten-fold increase in measurements, there is a reduction in eigenvalues by an order of magnitude.
The eigenvalues associated with the fifteen uninformed directions remain unaffected by $N$, appearing as a solid line in the plot.


\subsection{Higher-Dimensional Linear Gaussian Examples}\label{sec:high-dim-linear-example}

\begin{figure}[htbp]
\centering
  \includegraphics[width=0.5\linewidth]{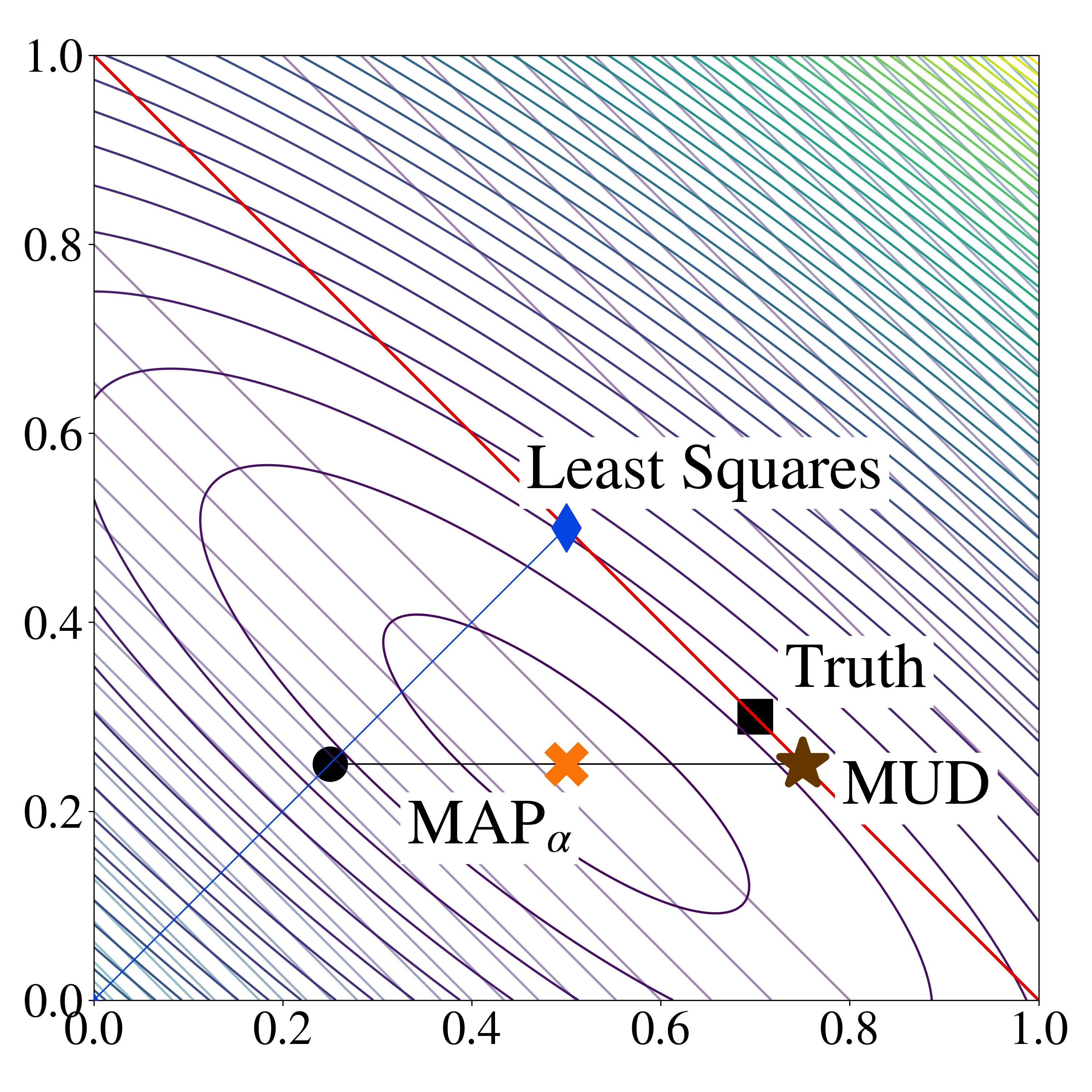}
\caption{
	A comparison of Least Squares, MAP, and MUD estimates relative to a true reference parameter along the solution contour defined by $Q^{-1}(\observedMean)$.
    Pictured are the contour lines involved in Tikhonov regularization defined in Table~\ref{tab:func_comparisons}. 
    While the MUD point will always lie on the solution contour, the MAP solution can be influenced by the strength of the regularization parameter $\alpha$ so that it lies along the line connecting the initial mean and the MUD point.
}
\label{fig:regularization-comparison}
\end{figure}

We first describe the relationship of MUD, MAP, and least squares estimates to the set-valued inverses of $Q$ in order to establish a conceptual framework for interpreting the numerical results that follow.
Throughout this discussion, we refer extensively to Figure~\ref{fig:regularization-comparison}, which builds upon Figure~\ref{fig:regularization}, to make these ideas more clear.
In Figure~\ref{fig:regularization-comparison}, a prior/initial covariance is chosen such that the MAP estimate is approximately halfway between the initial mean and the contour of $Q$ on which both the MUD and least squares estimates exist.
A reference parameter $\paramref = (0.75, 0.25)$, labeled as ``Truth'' is also shown on the contour.
Note that the MAP estimate is labeled with a subscript $\alpha$.
In these examples, $\alpha$ is a scalar multiple of the prior/initial covariance.
Thus, $\alpha$ is a hyper-parameter that determines, in a sense, the ``strength of prior beliefs'' as quantified by the eigenvalues in the prior/initial covariance. 

As illustrated in Figure~\ref{fig:regularization-comparison}, the MUD point always exists on the ``solution contour'' defined by $Q^{-1}(\observedMean)$ regardless of the initial covariance.
This is in fact guaranteed by Theorem~\ref{thm:MUD_existence_uniqueness}(b).
We use the term ``solution contour'' in this case because if $A$ is a $p$-to-$m$ full rank linear map with $m<p$, then $Q^{-1}(\observedMean)$ exists as a $(p-m)$-dimensional linear hyperplane in $\pspace$.
This means that the MUD point retains the ``predictive precision'' of a least squares solution to the inverse problem (i.e., a parameter that minimizes the data-mismatch term $\norm{\observedMean - Q(\param)}_{\observedCov^{-1}}^2$) while incorporating the flexibility of prior beliefs in directions not informed by the QoI.
This is in fact expected given the roles of $\initial$ listed in Section~\ref{subsec:defns} coupled with Theorem~\ref{thm:MUD_existence_uniqueness}.
For under-determined or ill-conditioned problems, this suggests that ``good'' prior beliefs may be used to produce a MUD point that is more accurate (i.e., closer to truth) than a least squares solution as illustrated in Figure~\ref{fig:regularization-comparison}.
This is further explored in the examples below involving high-dimensional linear maps.

By contrast, the MAP point exists on a line connecting the initial mean, $\param_0$, and the generalized contour defined by $Q^{-1}(\observedMean)$ as illustrated in Figure~\ref{fig:regularization-comparison}.
The direction of this line is discovered by substituting \eqref{eq:map_cov_analytical} into \eqref{eq:map-point-analytical}, which reveals that the direction of the line is orthogonal to the nullspace of the image of $A$ under $\initialCov$; i.e., $\mathcal{N}(\initialCov A)^\perp$.
In fact, this line intersects the generalized contour defined by $Q^{-1}(\observedMean)$ precisely at the MUD point.
If one parameterizes the line between $\param_0$ and $\param^{\text{MUD}}$, then one can also identify $\param^{\text{MAP}}$ as a convex sum of these two points.
The weights of this convex sum, which determine the position of the MAP point on this line, are determined by the ``precision of data'' (i.e., on $\observedCov$) and the ``strength of prior beliefs'' (i.e., on $\Sigma_{\text{init}}$).
This is apparent both in Figure~\ref{fig:regularization-comparison} and also by comparing the location of the MAP point in the top right plot of Figure~\ref{fig:regularization} to the line segment connecting the initial mean to the MUD point in the bottom right plot of this same figure.
The impact of this is also explored in the following examples where we tune an $\alpha$ hyper-parameter appearing as a multiplicative factor in $\Sigma_{\text{init}}$.

\begin{figure}[htbp]
  \includegraphics[width=0.475\linewidth]{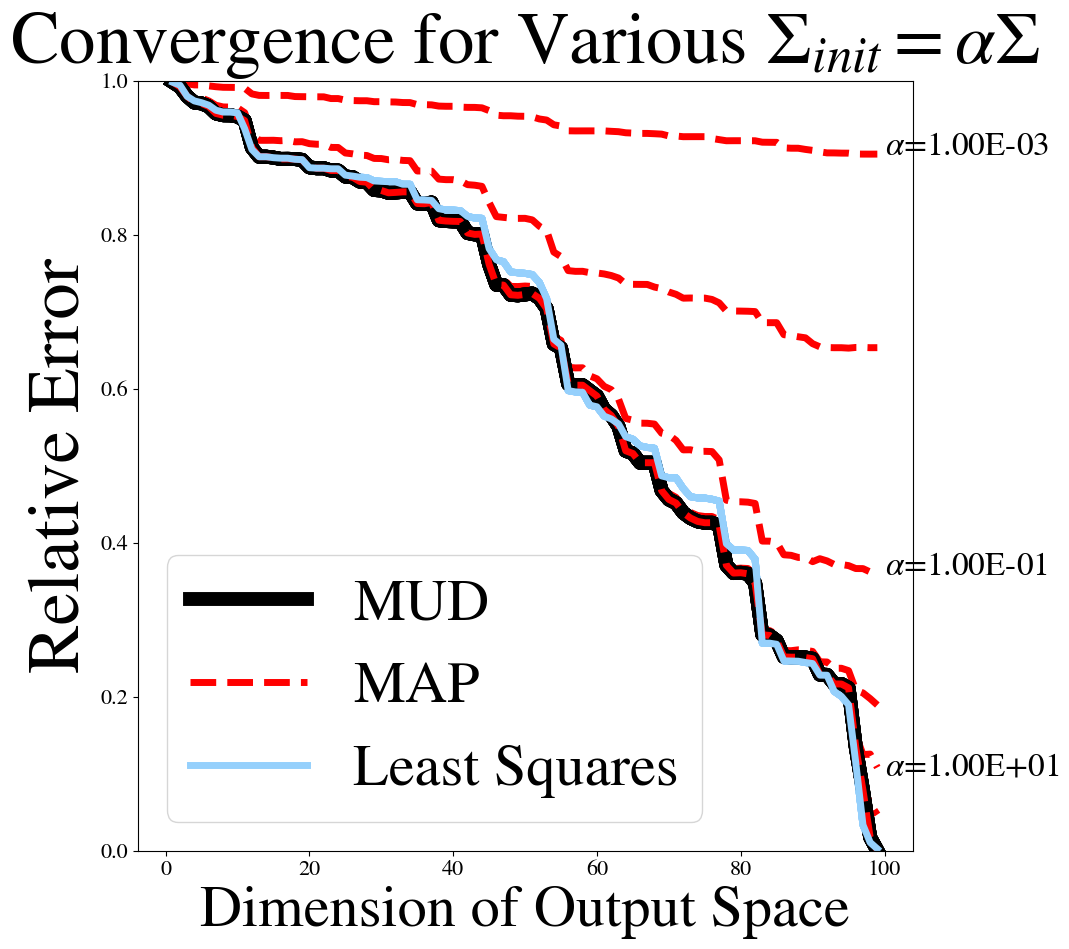}
  \includegraphics[width=0.475\linewidth]{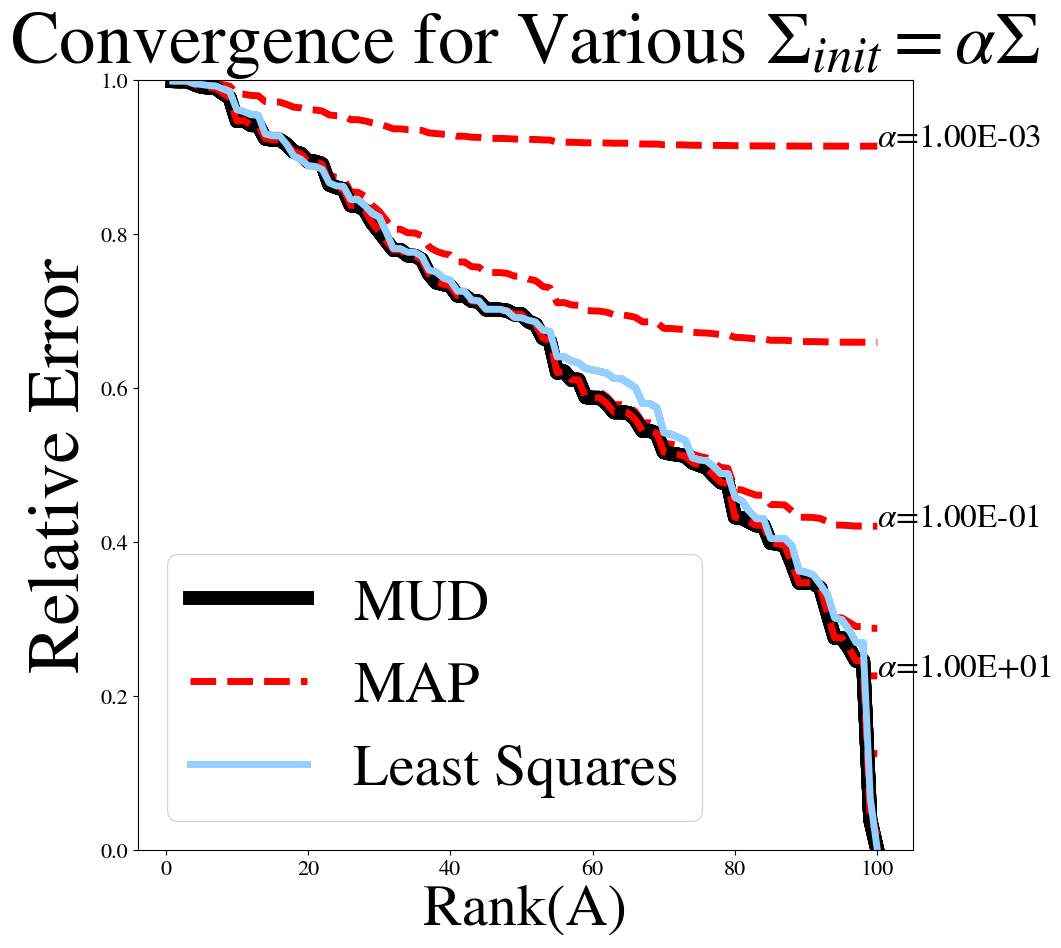}
\caption{
	Relative errors between $\paramref$ and (i) the least squares solution obtained through {\tt numpy}'s {\tt linalg.pinv} module, (ii) the closed-form solution for the MUD point given in Eq~\eqref{eq:mud-point-analytical-final}, and (iii) the MAP point.
  (Left): Error for increasing dimensions of $D$ for $A$ taken to be a Gaussian Random Map.
  (Right): Error for increasing row-rank of $A$, generated with Gaussian vectors and a SVD.
}
\label{fig:lin-error}
\end{figure}

\subsubsection{Impact of Output Dimension}


We consider QoI defined by $A\param+\mathbf{b}$ for $A\in\RR^{m\times 100}$ where $m=1,2,\ldots,100$ to demonstrate how the various estimates of a true parameter $\paramref$ are impacted by the number of available QoI.
To generate the matrices, we first generate $10,000$ i.i.d. random numbers from a $N(0,1)$ distribution that are arranged into a reference $\RR^{100\times 100}$ matrix.

The same distribution is also used for generating the components of the 100-dimensional vectors defining a reference bias vector $\mathbf{b}$ and reference parameter $\paramref$.
A multivariate Gaussian distribution is used for the initial density, with zero mean and $\initialCov$ chosen as a diagonal covariance with random entries drawn from $U[0.5, 1.5]$ and sorted in descending order. 
The prior density is also a zero mean multivariate Gaussian distribution.
However, to demonstrate the impact of the strength of prior beliefs on the MAP point, we choose the prior covariance to be $\alpha\initialCov$  for $\alpha=0.001, 0.01, 0.1, $ and $10$.
Here, smaller values of $\alpha$ correspond to a ``stronger'' belief in the prior since the prior density becomes more concentrated near the prior mean.

To study the impact of dimension on the MUD, MAP, and least squares estimates, we solve a
sequence of inverse problems by truncating the rows of the reference matrix and bias vector.
%
The results are summarized in the left plot of Figure~\ref{fig:lin-error}, which shows convergence towards $\paramref$ for all the problems considered with the exception of several MAP estimates corresponding to strongly-held beliefs in the prior.

We note that the MUD solution is the same for all choices of $\alpha$ and corresponds to the same level of accuracy that the MAP point achieves when $\alpha$ is chosen to be large.
In other words, the MUD point is not impacted by a scaling of the initial covariance, providing \emph{consistent} solutions which demonstrate levels of accuracy that MAP points only exhibit for larger values of scaling factors.

Of interest is also that the MUD point can sometimes out-perform the least squares estimate while generally achieving similar levels of accuracy.
This suggests that the MUD point has several favorable qualities.
Not only is it robust to the specification of prior assumptions, but it manages to offer the flexibility of incorporating good prior specifications without paying the additional cost of hyper-parameter optimization (i.e., choosing an appropriate $\alpha$) that would be required for the MAP estimates to achieve comparable results.

While omitted in the interest of space, if $\initialCov$ is chosen as $\alpha I$, where $I$ denotes the identity matrix of appropriate dimensions, then the MUD point will always agree with the least squares estimate.
Taking these results together, this implies that only a good ``relative spatial structure'' of prior beliefs is required to improve the MUD point's accuracy over both MAP and least squares estimates.

\subsubsection{Impact of Rank: One Hundred (Deficient) $\RR^{100 \times 100}$ matrices}

Here, we investigate whether the previous dimension-dependent example results also apply to matrices $A$ which are of a fixed dimension but possess varying rank.
This is of interest in applications where many QoI are available to construct an operator but a great deal of redundancy may be present in the data collected, and feature-engineering new quantities is somehow prohibitive (perhaps due to gradient estimation).


The rank of $A$ corresponds to the number of unique directions of information present in the operator, i.e., how many directions in the parameter space are informed by the QoI map.
The operators in the previous example were all full rank, so the dimension of each map also corresponded to the rank of $A$.
When $A$ is rank-deficient, $\predictedCov$ is non-invertible, so we must modify the form of \eqref{eq:mud-point-analytical-final} to substitute a pseudo-inverse for the predicted covariance.

In this example, the dimension of the data space remains fixed at $m=100$ across all experiments.
However, we sequentially increase the row-rank of $A$ from $r=1, \ldots, 100$.
To control the rank of $A$, we first construct a reference $\RR^{100\times 100}$ matrix as in the previous example using 10,000 i.i.d. $N(0,1)$ random numbers.
We then compute a singular value decomposition of this reference matrix of the form $USV^\top$ and construct 100 rank-1 matrices of the form $A_i=\mathbf{u}_i s_i \mathbf{v}_i^\top$ for $i=1,\ldots,100$ where $\mathbf{u}_i$ and $\mathbf{v}_i$ denote the $i$th columns of $U$ and $V$, respectively and $s_i$ denotes the $i$th singular value.
Then, we analyze the impact of $A = \sum_i^r A_i$ for $r=1,\ldots,100$.
Aside from the differing construction of $A$, the rest of the choices involved in the experiment ($\paramref$, the reference bias vector, and the distributions involved) is identical to the previous example.

In the right plot of Figure~\ref{fig:lin-error}, we
again find that the MUD point is generally as accurate as the least squares estimate, but incorporates an initial description of uncertainty, which may allow it to outperform the least squares estimate.
Also, we again see that the MAP estimates are impacted by the strength of prior beliefs.

\section{Data-constructed QoI maps II: Data Clouds and Principal Component Analysis}\label{sec:PCA}

In section \ref{sec:data-maps} the $Q_{WME}$ is introduced as a way to incorporate data from repeated measurements and reduce the variance in the MUD point estimates.
Here, we address the case of potentially non-repeated measurements taken from a system that varies over a temporal and/or spatial domain.

We utilize the ubiquitous Principal Component Analysis (PCA).
Introduced first by \citep{pearson1901liii,hotelling1933analysis}, PCA is a way to reduce dimensionality of a large set of correlated data by transforming the data into a new set of variables known as principal components, which are uncorrelated and when ordered, contain the maximum amount of variation in the data per component.
For more comprehensive reviews of the PCA method, we refer the interested reader to \citep{shlens2014tutorial, jolliffe2002principal}.
The PCA transformation can be written in the form

\begin{linenomath*}
\begin{equation}\label{eq:pca}
	Y = XP,\quad X, Y \in \mathbb{R}^{s\times n},\quad P \in \mathbb{R}^{n\times n}.
\end{equation}
\end{linenomath*}

Here, the data matrix, $X$, contains $s$ samples of $n$ data points each.
Using a change of basis matrix $P$, $X$ is transformed to a new matrix $Y$. We summarize below some of the main conclusions of the PCA transformation.

\begin{itemize}
	\item $P$ defines a linear transformation to an orthonormal basis, given by the columns $p^{(\ell)}\in \mathbb{R}^n, 1 \leq \ell \leq n$ of the matrix $P$ that define the principal components of $X$.
	\item $P$ diagonalizes the covariance matrix $C_{X} = XX^T$. 
	\item The matrix $P$ and the popular Singular Value Decomposition (SVD) transformation are intimately linked. Namely, for the SVD decomposition of $X = U{\Sigma}V^T$, the columns of $V$ are the principal components of $X$.
\end{itemize}

Returning to the scenario proposed in section~\ref{sec:data-maps}, suppose the $\mathcal{M}_j$ measurement devices now collect data over space and time, with each taking $N_j$ measurements as before.
Define any ordering of these $n = \sum_{j=1}^m N_j$ data points $\{z_i\}_{i=1}^n$ (the ordering can be arbitrary as it does not impact the results obtained via PCA).
Similar to equation (\ref{eq:obs_data_error}), we now let $d_i$ equal the $i$th measurement datum, which is assumed polluted by i.i.d. additive Gaussian errors from a $\mathcal{N}(0,\sigma_i)$ distribution.
Furthermore, let $\mathcal{M}_{k,i} = \mathcal{M}(\lambda_k;z_i)$ be the $i$th measurement for the $k$th simulated sample.
Assuming $s$ samples are collected, the matrix $X \in \mathbb{R}^{s\times n}$ of Z-scored residuals for a sample set is now defined component-wise as

\begin{linenomath*}
\begin{equation}\label{eq:matrix_res}
	X_{k,i} = \frac{\mathcal{M}_{k,i} - d_i}{\sigma_i}.
\end{equation}
\end{linenomath*}

We define $Q_{PCA}$ component-wise as

\begin{linenomath*}
\begin{equation} 
	(Q_{PCA})_\ell(\lambda_k) = \sum_{i=1}^n p^{(\ell)}_i \frac{\mathcal{M}(\lambda_k; z_i) - d_i}{\sigma_i}, \quad 1 \leq \ell \leq n,
\end{equation}\label{eq:q_pca}
\end{linenomath*}
where $p^{(\ell)}$ is the $\ell$th principle component of $X$.

Similar to $Q_{WME}$, the $Q_{PCA}$ map is computing a weighted average of residuals.
However, the weighting of these residuals is according to the new set of basis vectors defined by the principal components of $X$.
Since each component of $Q_{PCA}$ is still a normalized (in the 2-norm) combination of the Z-scored residuals, the observed distribution $\observed$ is still given by an $\mathcal{N}(0,1)$ distribution as with the $Q_{WME}$ map.
Note that the choice of ordering of the measurements in the data matrix $X$ is irrelevant, since the PCA does not depend on column order.
Finally, observe we can define up to $n$ components to our $Q_{PCA}$ map.
However, in practice, we only take up to the first $m$ components that capture a user-specified percentage of variance in the original data set $X$. 
In fact, if the data are sensitive to all parameters present in our inverse problem, we expect the number of components, $m$, to be equal to the dimension of our paremter space, $p$. 
However, we see in the following examples situations where this may not be the case.
In these examples, we turn to the diagnostic $\mathbb{E}(r)$ as an important measure of the quality of the updated density and thus the reliability of $\mudpt$.

\section{Spatial and Temporal Data Examples}\label{sec:PCA_examples}

We now use MUD points as parameter estimates for PCA constructed maps using simulated noisy temporal and spatial data associated with solutions to differential equations.
The previously derived closed form expressions for $\mudpt$ do not apply in these examples since the measurement maps from parameters to data are nonlinear.
Instead, in each example we use a fixed set of i.i.d.~samples drawn from the initial density to approximate the updated density and subsequently choose the sample that maximizes this approximation.
In each case, we monitor $\mathbb{E}(r)$ to check for the validity of the predictability assumption.

\subsection{Spatial Data Example: Poisson's Equation with Uncertain Boundary Condition}

We first consider the aggregation of data over a spatial domain. In this problem, the uncertain model parameter is described by an uncertain parametrized function defining the boundary data to a stationary PDE, given by the Poisson problem:
\begin{linenomath*}
\begin{equation}\label{eq:pde-equation}
\begin{cases}
\hfill -\nabla \cdot \nabla u &= f(x), \quad\text{on } x\in \Omega, \\
\hfill u &= 0, \quad\text{ on } \Gamma_T \cup \Gamma_B, \\
\hfill \frac{\partial u}{\partial \mathbf{n}} &= g(x_2), \quad\text{ on } \Gamma_L, \\
\hfill \frac{\partial u}{\partial \mathbf{n}} &= 0, \quad\text{ on } \Gamma_R,
\end{cases}
\end{equation}
\end{linenomath*}
where $x=(x_1, x_2) \in \Omega = (0,1)^2$ is the spatial domain; $\Gamma_T$, $\Gamma_B$, $\Gamma_L$, and $\Gamma_R$, denote the top, bottom, left, and right boundaries of this domain, respectively, and $\frac{\partial u}{\partial \mathbf{n}}$ denotes the usual outward normal derivative.
The forcing function $f$ is taken to be $10\exp\left ( \norm{x - 0.5}^2 / 0.02 \right )$.

\begin{figure}[htbp] 
\centering
\includegraphics[width=0.95\linewidth]{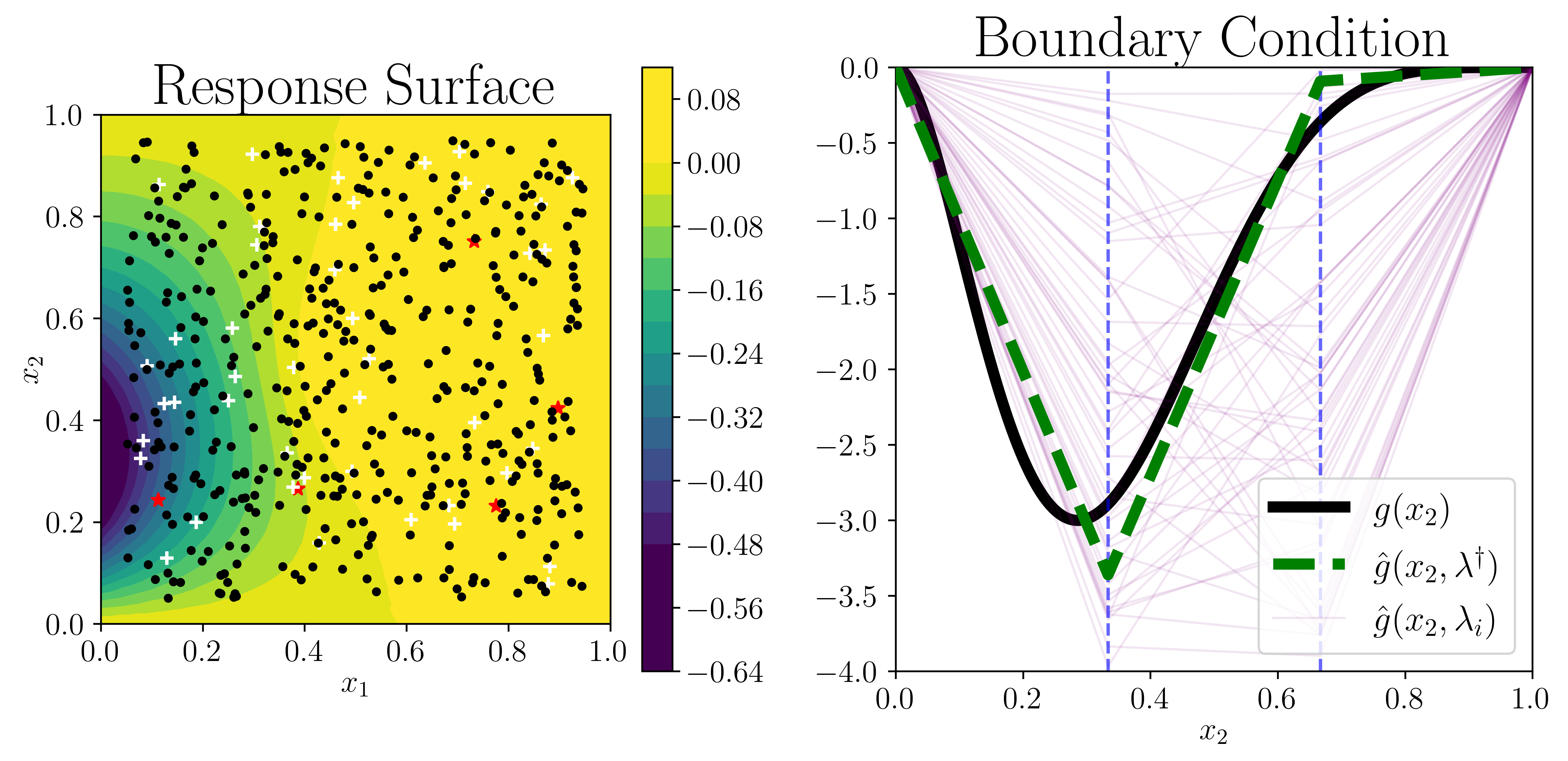}
\caption{(Left) Noiseless response surface to solving equation (\ref{eq:pde-equation}) using true $g(x_2)$ with locations of sensors used for collecting data. First five sensors used in red in stars, next set of 45 in white crosses, and the rest (450) are in black circles. (Right) Left boundary condition. True $g(x_2)$ (black) along with piecewise-linear spline approximations $\hat{g}(x_2,\lambda)$. The closest (in the 2-norm) possible spline to the ``true" boundary data is indicated by $\lambda^\dagger$ (dotted green). 50 different $\lambda_i$ samples from the initial uniform distribution are shown in faded purple).
}
\label{fig:pde_surface}
\end{figure}

The goal is to use noisy data to estimate the boundary data $g(x_2)$.
To generate the noisy data, a ``true" $g(x_2)\propto x_2^2(x_2-1)^5$ is constructed, with a constant of proportionality chosen to produce a minimum of $-3$ at $x_2=\frac{2}{7}$.
Equation (\ref{eq:pde-equation}) is solved using piecewise-linear finite elements on a triangulation of a $36\times36$ mesh.
Random noise is then added to every degree of freedom of this reference solution, and the spatial data are subsequently computed from a fixed set of 5, 50, and 500 randomly placed sensors in the subdomain $(0.05, 0.95)^2 \subset \Omega$ (Figure ~\ref{fig:pde_surface}, left).

To construct a finite-dimensional parameter space describing the initial uncertainty of $g(x_2)$, piecewise-linear continuous splines, $\hat{g}(x_2)$, are used to approximate $g(x_2)$.
The locations of the first and last knots are fixed at the endpoints of the boundary with values assumed to be 0, and furthermore it is assumed that $g$ is non-positive and bounded below by $-4$.
Thus, the uncertainty is described by the values of the splines at the interior knot points chosen as the equispaced points $1/3$ and $2/3$.
This defines a finite-dimensional parameter space described by $\pspace = [-4,0]^2$.
We generate $1000$ samples from an initial uniform density on $\pspace$ to (1) generate random spline functions and compute the (noise-free) data from solutions associated with these splines; and (2) estimate the push-forward and updated densities along with the MUD estimate of the boundary data, $\hat{g}(x_2, \lambda^{MUD})$.
Note that since we approximate the ``true" boundary data with linear splines, the MUD estimate will {\em never} exactly equal the noiseless boundary data.
To measure the accuracy of a given estimate then, we compare it with the closest (in the 2-norm) possible spline to the ``true" boundary data.
We refer to the values of $\lambda$ that produce this value as $\lambda^{\dagger}$.

\begin{figure}[htbp] 
	\centering
	\includegraphics[width=\textwidth]{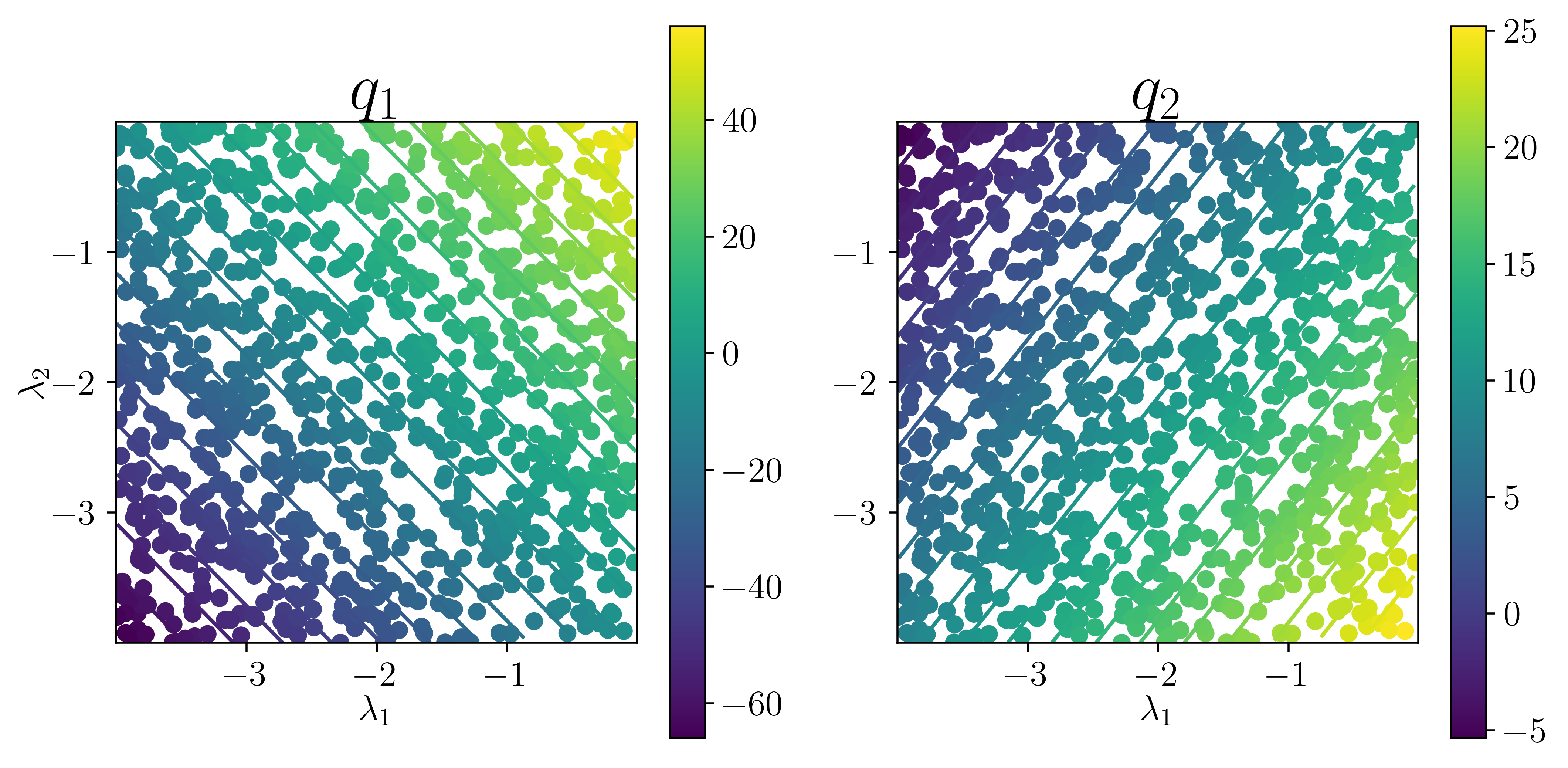}
	\caption{Quantities of interest learned by the $Q_{PCA}$ map when using all 500 sensors for the first (left) and second (right) principal component. Note how the contours show that the $Q_{PCA}$ map builds orthogonal components that lead to a well posed inversion problem.}
	\label{fig:learned_qoi}
	\centering
\end{figure}

For each sample, the $Q_{PCA}$ map in equation (\ref{eq:q_pca}) is constructed using two principal components.
The value of the learned QoI map for each parameter sample is seen in Figure \ref{fig:learned_qoi}.
Note the orthogonal contour structure of the two components demonstrates how the $Q_{PCA}$ map aggregates data from a high-dimensional output space into an essentially linear bijective map with a Jacobian that is well-conditioned for inversion.

\begin{figure}[htbp] 
\centering
   \vspace{-5mm}
   \includegraphics[width=0.82\linewidth]{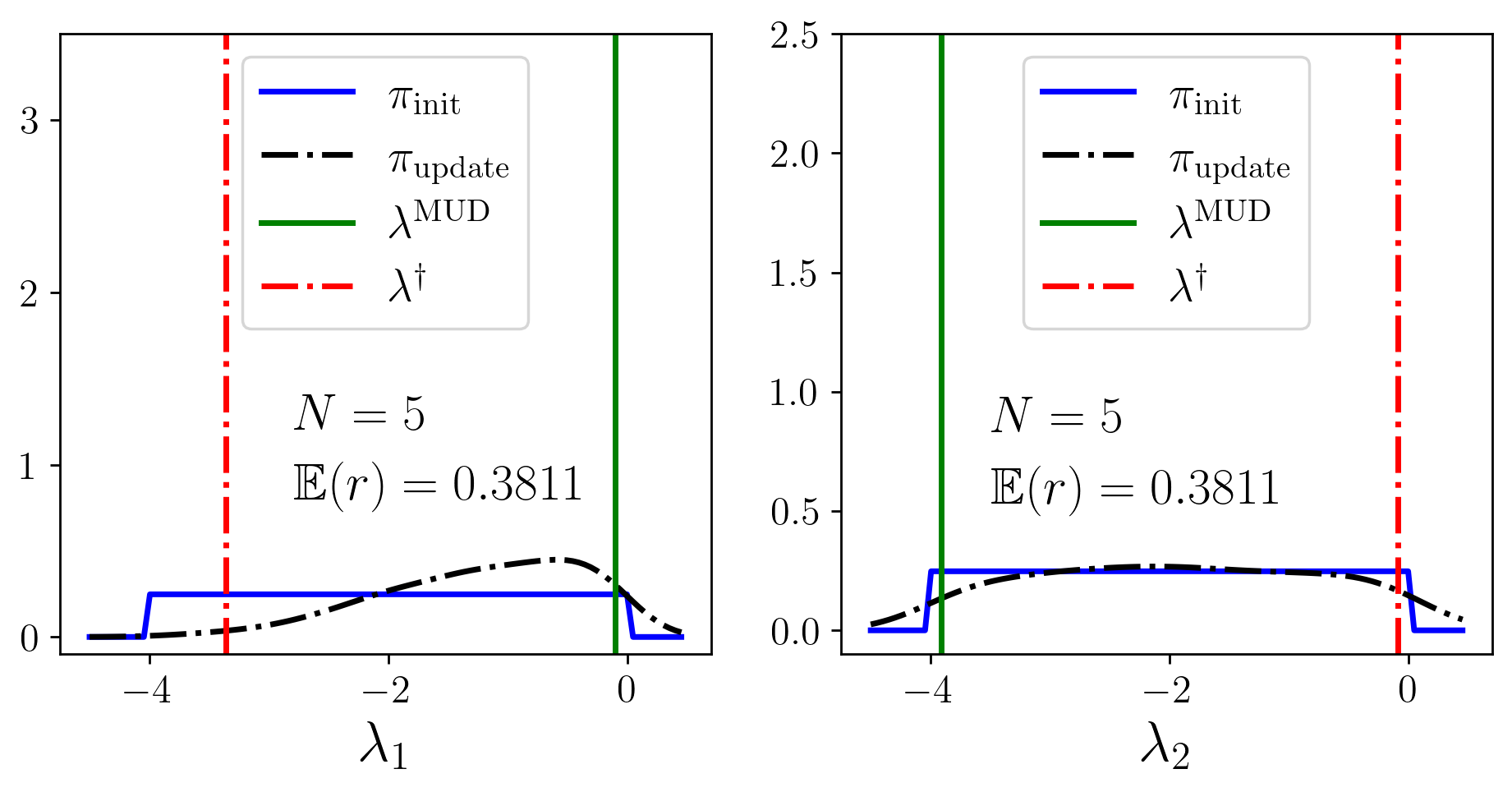}
   \includegraphics[width=0.82\linewidth]{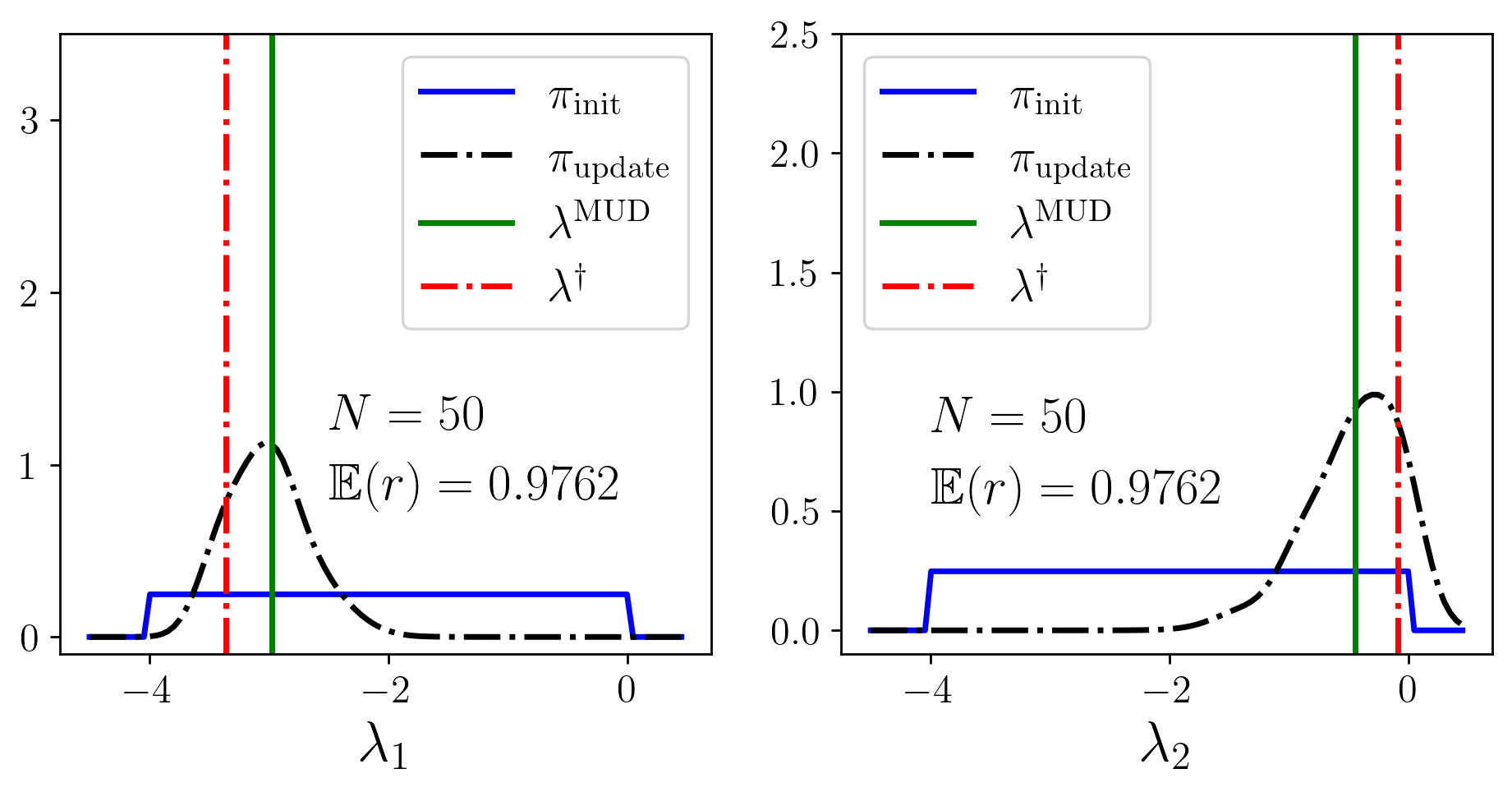}
   \includegraphics[width=0.82\linewidth]{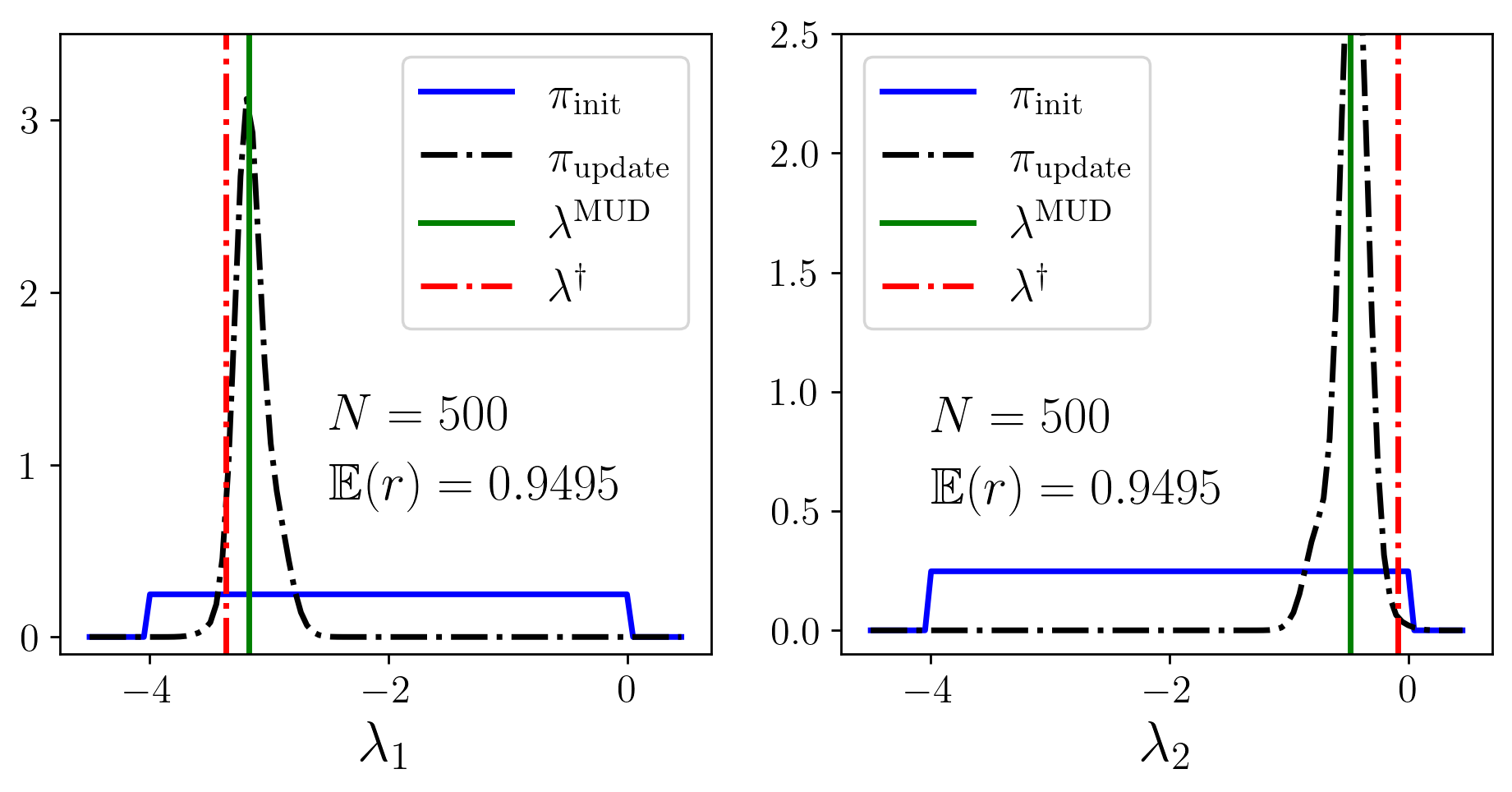}
   \caption{Initial (blue) and updated (dotted-dashed, black) distributions for $\lambda_1$ (left plots) and $\lambda_2$ (right plots) using $N=5 \text{ (top), }50 \text{ (middle), }500 \text{ (bottom)}$ sensors. Expectation value of the ratio of the observed to the predicted, which should be $\approx 1$, is shown as well per case.}
   \label{fig:pde2d-comparison-convergence}
\end{figure}

The plots in Figure ~\ref{fig:pde2d-comparison-convergence} show how the quality of MUD estimates is improved as more sensor data are included.
Note how for only the first five sensors (red stars in Figure ~\ref{fig:pde_surface}, left), the $Q_{PCA}$ map does not produce a good estimate for either parameter.
In fact, practically no update is made in $\lambda_2$, while a small update is made in $\lambda_1$, but in the incorrect direction.
The location of the first five sensors explains this behavior as there are only two sensors near the location of the knot controlled by $\lambda_1$ on the boundary, while no sensors are near the location on the boundary where the $\lambda_2$ knot location is.
It is important to stress the importance of the diagnostic $\mathbb{E}(r)$ here.
If we had no knowledge of the ``true" boundary data, and no $\lambda^\dagger$ was readily available for comparison,  the results of using five sensors may lead us to believe that a meaningful update was made for $\lambda_1$.
However, the diagnostic of $\mathbb{E}(r) = 0.3811$ tells a different story, indicating that the predictability assumption is likely being violated in the construction of $\updated$.
Once $N=50$ sensors are used (red stars and white crosses in Figure ~\ref{fig:pde_surface}, left), the diagnostic value jumps to within an acceptable range of $\mathbb{E}(r) = 0.9762 \approx 1$, and the $Q_{PCA}$ map is able to update initial beliefs in both parameters well.
Using $N=500$ sensors leads to a reduction in the variance of the MUD point estimate as distributions for both parameter values peak around a sample that is very close to $\paramref$.

\subsection{Temporal Data Example: ADCIRC with Uncertain Wind Drag}\label{sec:adcirc}

For the second example, the aggregation of data over time into a $Q_{PCA}$ map is considered.
In this problem, the uncertain parameter relates to a coefficient that determines how unresolvable dynamics of a PDE system are modeled.
We begin with the Shallow Water Equations (SWE), which are a depth-averaged approximation to the Navier-Stokes equations commonly used in coastal circulation and flooding modeling to predict peak storm-surge due to extreme weather events \citep{vreugdenhil1994numericala}. 
They can be expressed as 
\begin{linenomath*}
\begin{equation}\label{eq:swe}
	\begin{cases}
\hfill \frac{\partial \zeta}{\partial t} + \nabla \cdot (\mathbf{U} H) &= 0, \\
\hfill \frac{\partial \mathbf{U}}{\partial t} + \mathbf{U} \cdot \nabla \mathbf{U} + f \mathbf{k} \times \mathbf{U} &= -\nabla\left[ \frac{p_s}{\rho_0} + g\zeta\right] + \frac{\boldsymbol\tau_s - \boldsymbol\tau_b}{\rho_0 H},
	\end{cases}
\end{equation}
\end{linenomath*}
for unknown free surface elevation, $\zeta=\zeta(x, y, t)$ with respect to mean sea level and depth-average velocity (averaged over the height of the water column), $\mathbf{U} = \mathbf{U}(x, y, t)$. 
In equation (\ref{eq:swe}), $H$ is the height of the water column, $f$ the Coriolis parameter, $p_s$ the atmospheric pressure, $\rho_0$ the reference density of water, $g$ the gravitational constant, $\boldsymbol\tau_s$ the surface stress, and $\boldsymbol\tau_b$ the bottom stress.
The surface stress,
\begin{linenomath*}
\begin{equation}\label{eq:surface_stress}
\boldsymbol\tau_s = \rho_s C_d \mathbf{u} ||\mathbf{u}||,
\end{equation}
\end{linenomath*}
includes the wind drag $C_d$, which is an effective (homogenized) parameter governing the transfer of momentum from winds to the water column, one of the primary drivers of storm-surge.
There are various forms for $C_d$ that depend on the physical properties of the system being modeled (e.g., the type of storm, presence of ice, etc.).
In this work, a popular generalization of Garratt's formula for $C_d$ is used \citep{garratt1977review}, whereby $C_d$ increases linearly with wind speed $u := ||\mathbf{u}||$ according to the uncertain parameter $\lambda_1$ and is ``cut off" when exceeding a threshold $\lambda_2$:

\begin{linenomath*}
\begin{equation}\label{eq:garratt}
C_d = \min\left[10^{-3} (.75 + \lambda_1 u), \lambda_2\right].
\end{equation}
\end{linenomath*}

Using time series of recorded water surface elevations, the goal is to determine the values of the wind drag parameters $(\lambda_1, \lambda_2)$.

\begin{figure}[htb] 
\centering
   \includegraphics[width=0.75\linewidth]{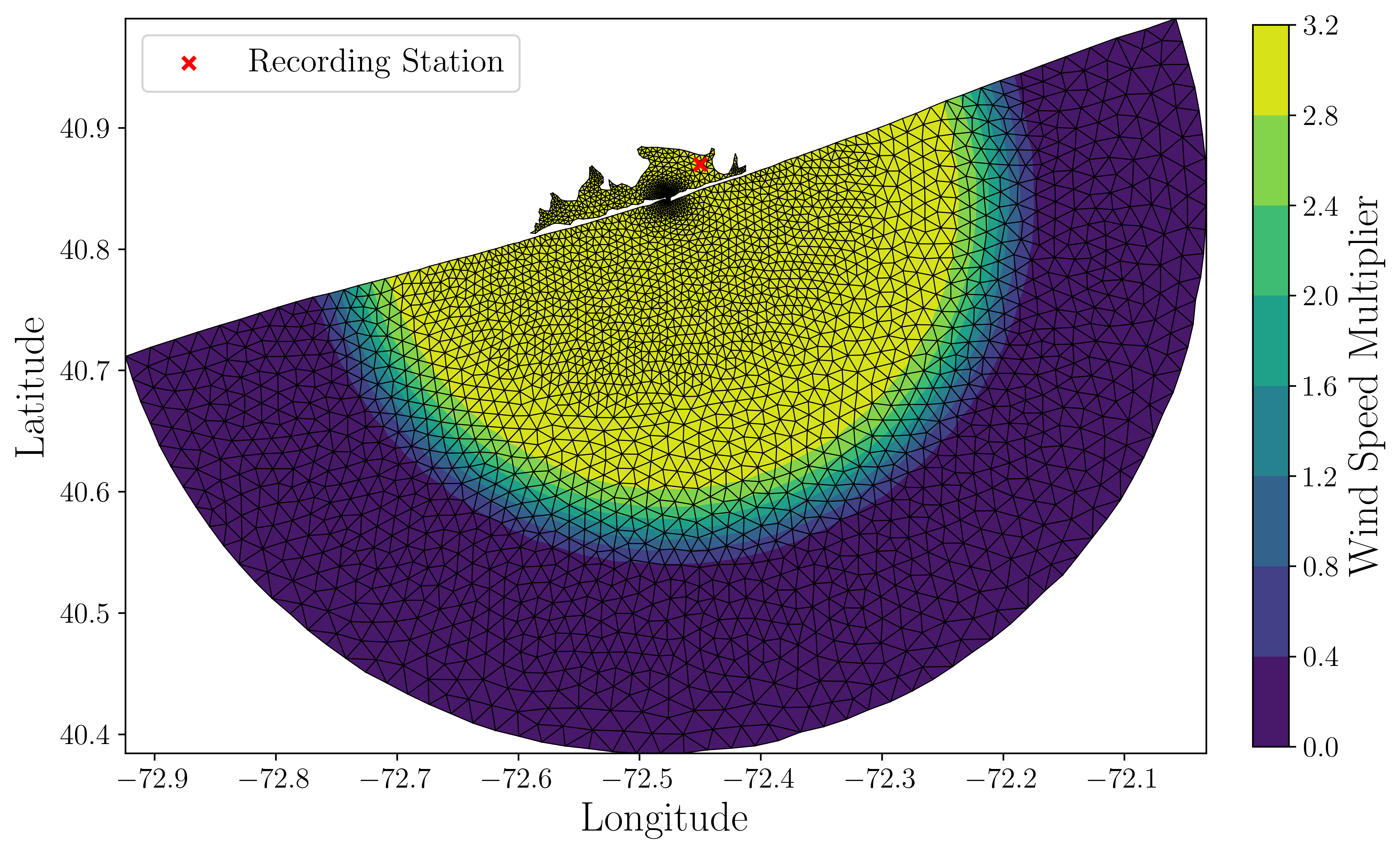}
   \includegraphics[width=0.85\linewidth]{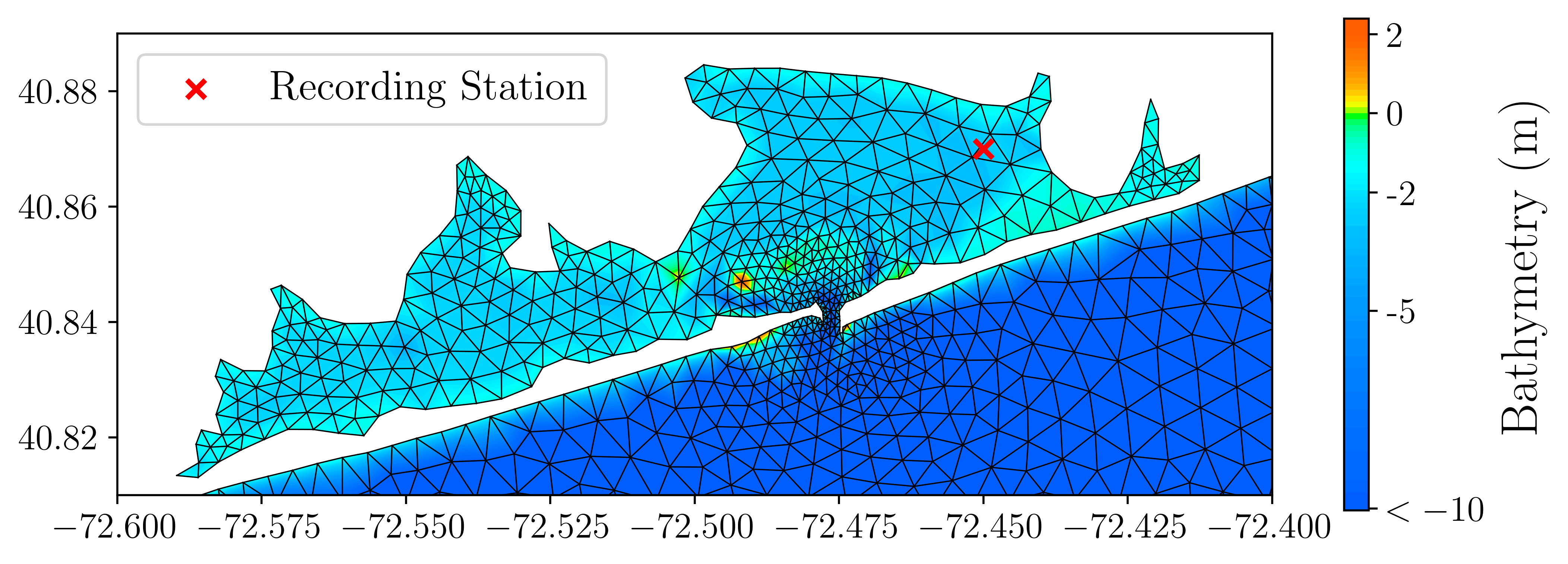}
   \caption{(Top) Shinnecock Inlet Mesh containing 5780 triangular elements. Contours show value of wind multiplier applied to scale winds up artificially near the inlet, tapering them off to zero at the boundaries.(Bottom) Bathymetry of inlet.}
 \label{fig:si-inlet}
\end{figure}

The ADvanced CIRCulation (ADCIRC) coastal ocean model is used to solve the SWEs (\ref{eq:swe}) \citep{luettich1992adcirc, westerink1992tide}.
ADCIRC uses a finite-element model of the SWE in which the Genralized Wave Continuity Equation \citep{lynch1979wave} is descretized in space using piecwise-linear elements on unstructured (triangular) grids.
The model is used for coastal engineering applications such as hurricane storm surge forecasting \citep{dietrich2013realtime}, hindcasting \citep{bunya2010highresolution, dietrich2010highresolution, dietrich2011hurricane} and uncertainty quantification \citep{butler2015definition, graham2015adaptive, graham2017measuretheoretic}. 

In this study, ADCIRC is configured to run using a test mesh based on the Shinnecock Inlet on the Outer Barrier of Long Island, NY, USA. 
External forcing for the model is given by tidal forcing reconstructed from the TPXO9.1 harmonic tidal constituents \citep{egbert2002efficient} using OceanMesh2D \citep{roberts2019oceanmesh2d}, constant air pressure of 1013 millibars, and free surface stress from winds computed from a $0.25\deg$ hourly CFSv2 10-m wind fields \citep{saha2014ncep} for a period of 16 days (29 December 2017 - 31 January 2018).
Winds are modified for the purposes of the numerical experiment to simulate a more extreme (Category 4) event, with winds scaled radially down to zero from the point of interest, i.e. the center of the inlet (see Figure ~\ref{fig:si-inlet}).

\begin{figure}[htbp] 
\centering
   \includegraphics[width=0.95\linewidth]{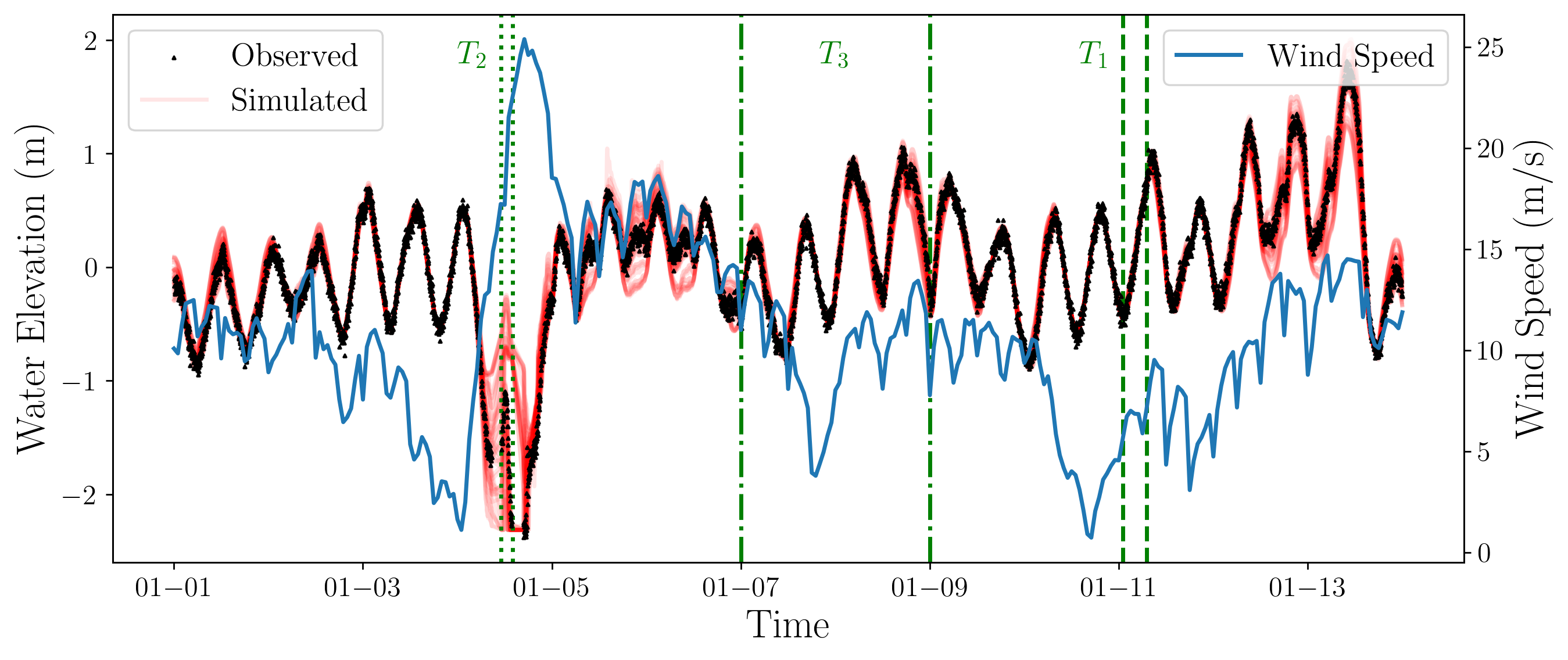}
   \caption{Time series of water elevations (left axis) for ``observed"  data (black triangles) and simulated data (faded red lines) along with wind forcing (right axis, blue line). The three time windows of data used are indicated in the vertical dashed ($T_1$), dotted ($T_2$), and dashed-dotted ($T_3$) green lines.}
   \label{fig:adcirc-res}
\end{figure}

To frame the PIP, first we assume that the uncertain parameters $(\lambda_1, \lambda_2)$ lie within $\pm 50\%$ of commonly used default values of $(0.067, 0.0025)$ \citep{letchford2009wind}.
This defines a finite-dimensional parameter space
\begin{linenomath*}
\begin{equation} 
	\pspace = [0.0335, 0.1105]\times [0.00125,0.00375]\subset\mathbb{R}^2. \label{eq:adcirc_lambda_domain}
\end{equation}
\end{linenomath*}
1000 samples are generated from a uniform distribution over $\pspace$ and pushed through our forward model, ADCIRC (see ~\ref{app:A} for more details on how to obtain simulation data).
Water elevation at an artificial recording station inside the inlet was recorded over a period of 14 days (1 January 2018 - 14 January 2018) at 3 hour intervals for each sample.
Since no real station data are available over the test mesh domain, we create observations by picking (and removing) a sample closest to the default value of $(0.067, 0.0025)$, and populate each measurement of this sample with i.i.d. $N(0, \sigma^2)$ noise, using $\sigma = 0.05$. 
The different time-series for the water elevation at the artifical recording station are shown in Figure ~\ref{fig:adcirc-res}.

\begin{figure}[htb] 
\centering
   \includegraphics[width=0.80\linewidth]{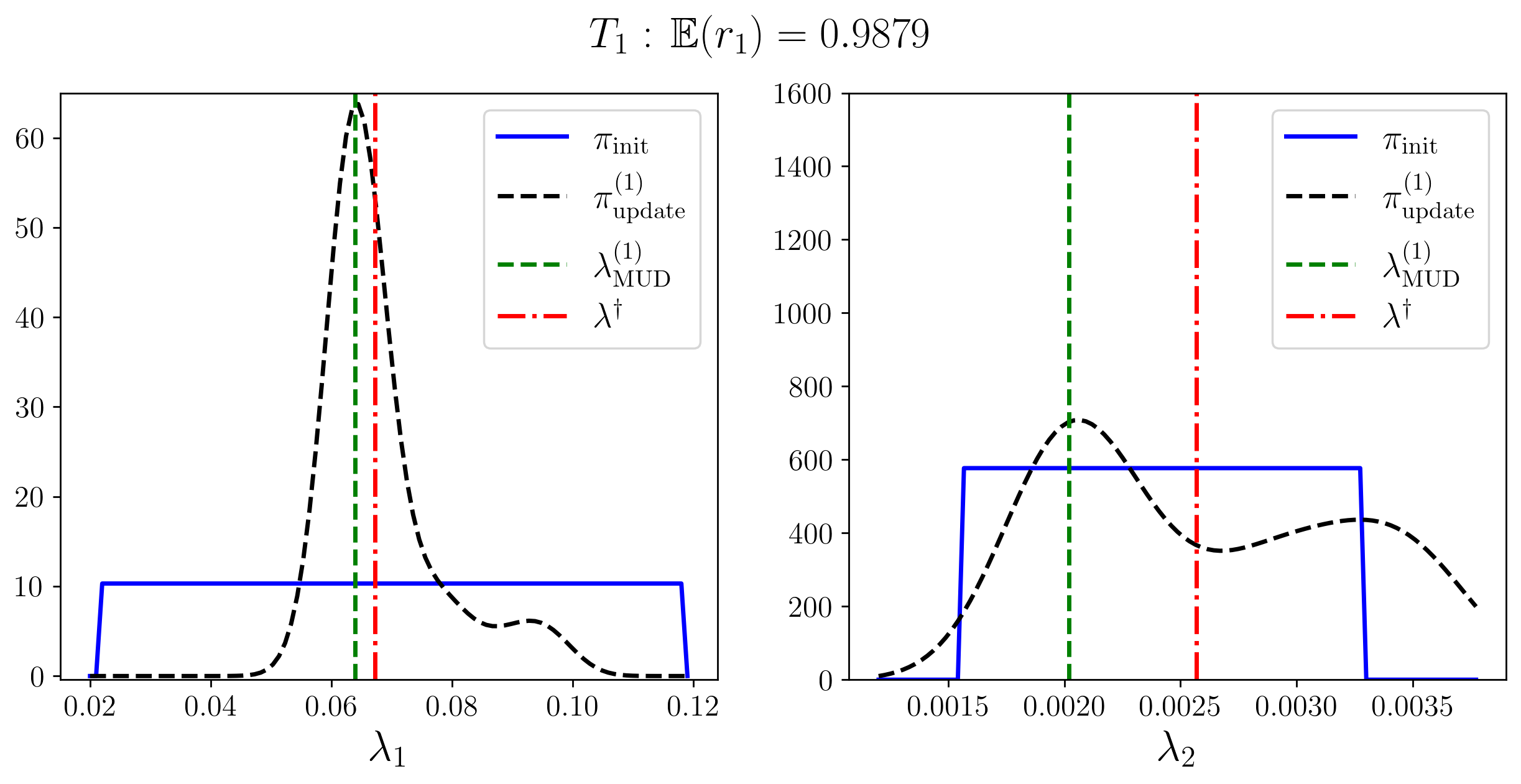}
   \includegraphics[width=0.80\linewidth]{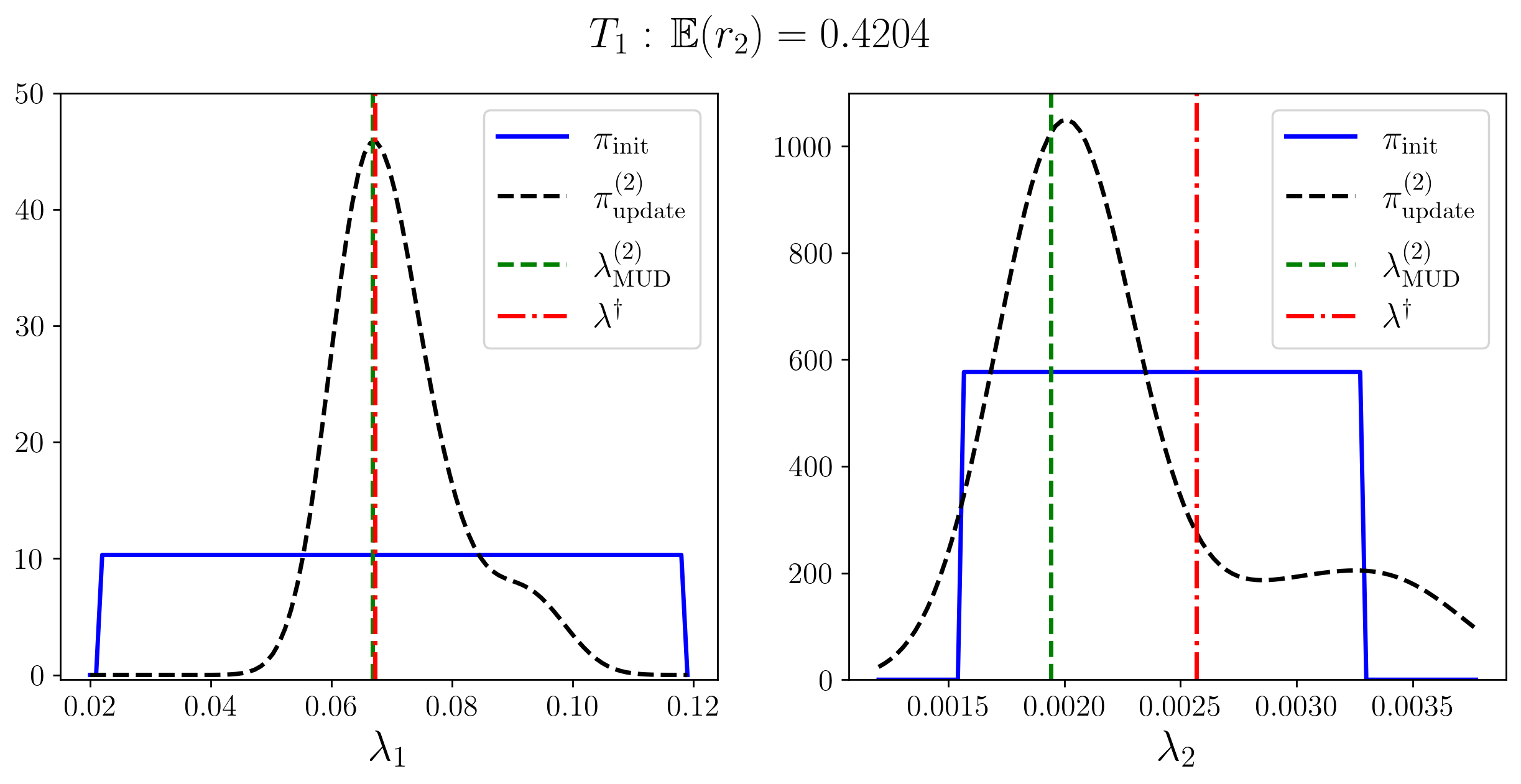}
   \caption{Updated distributions (dotted black line) and mud estimates (dotted green line) using a one component (top) vs two component (bottom) $Q_{PCA}$ map for $T_1$ (see Figure~\ref{fig:adcirc-res}) with $N = 119$ data points from January 11 01:00:00-07:00:00.}
   \label{fig:adcirc-t1}
\end{figure}

\begin{figure}[htb] 
\centering
   \includegraphics[width=0.80\linewidth]{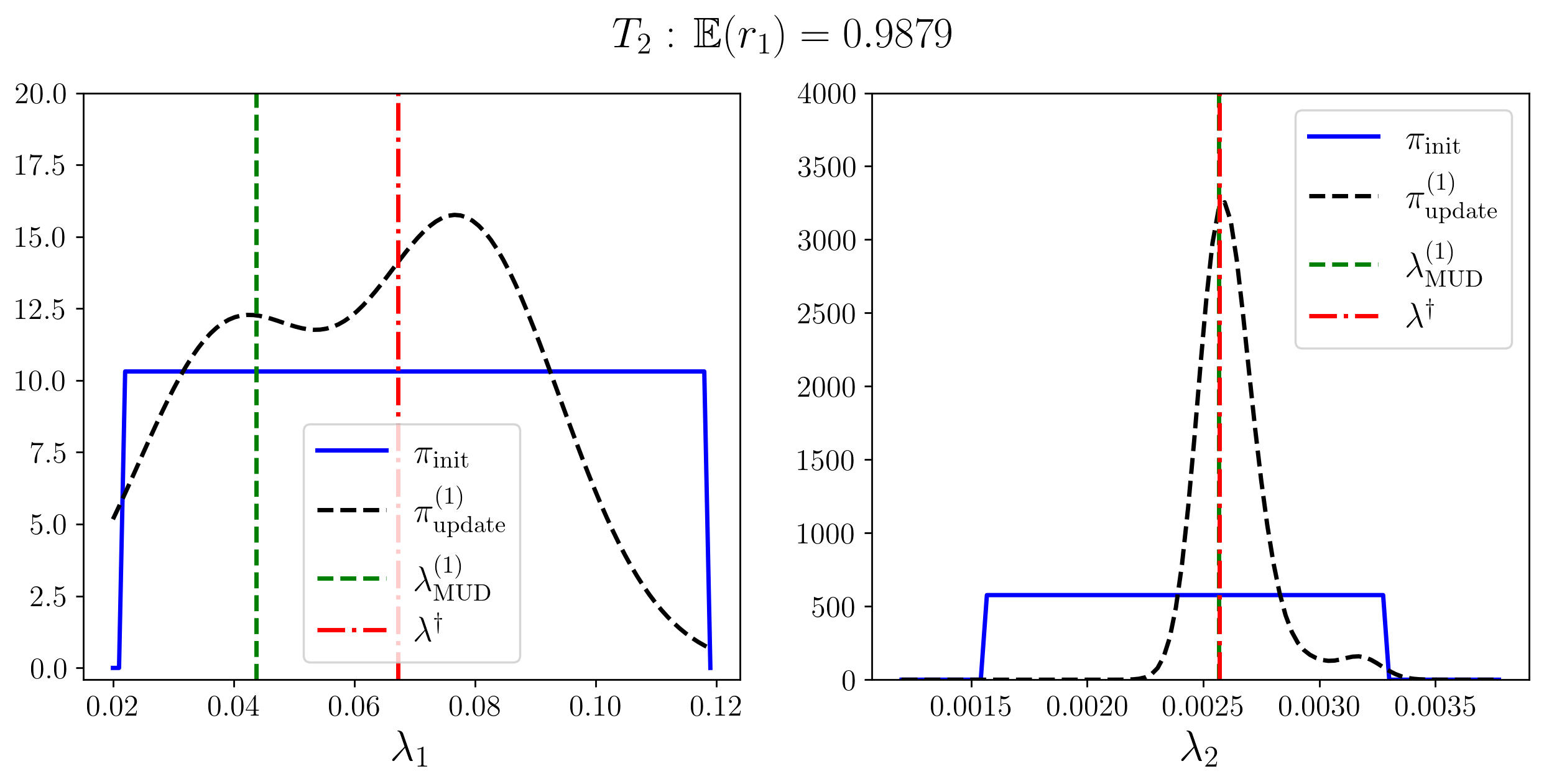}
   \includegraphics[width=0.80\linewidth]{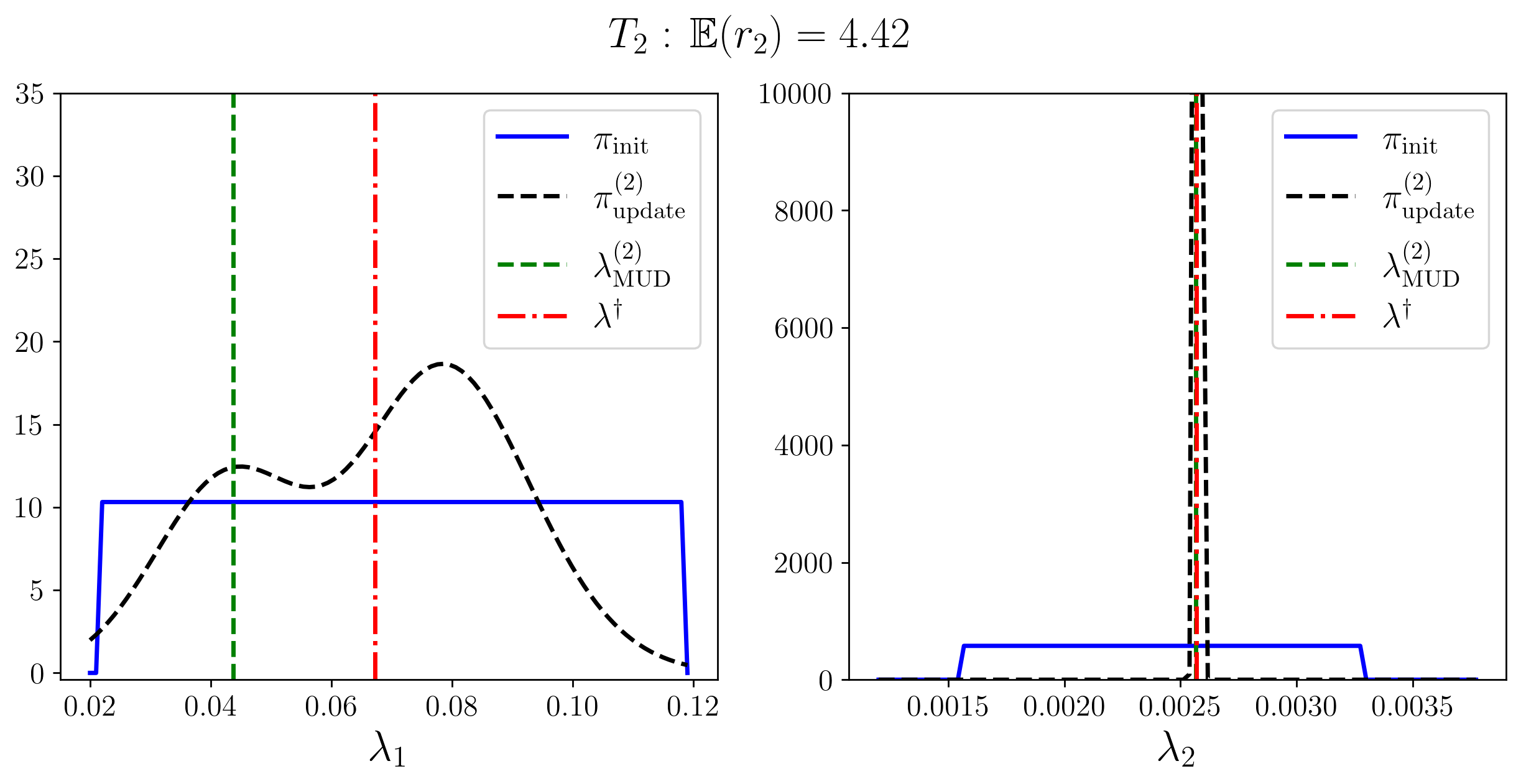}
   \caption{Updated distributions (dotted black line) and mud estimates (dotted green line) using a one component (top) vs two component (bottom) $Q_{PCA}$ map for $T_2$ (see Figure~\ref{fig:adcirc-res}) with $N = 53$ data points from January 04 11:00:00-14:00:00.}
   \label{fig:adcirc-t2}
\end{figure}

\begin{figure}[htb] 
\centering
   \includegraphics[width=0.80\linewidth]{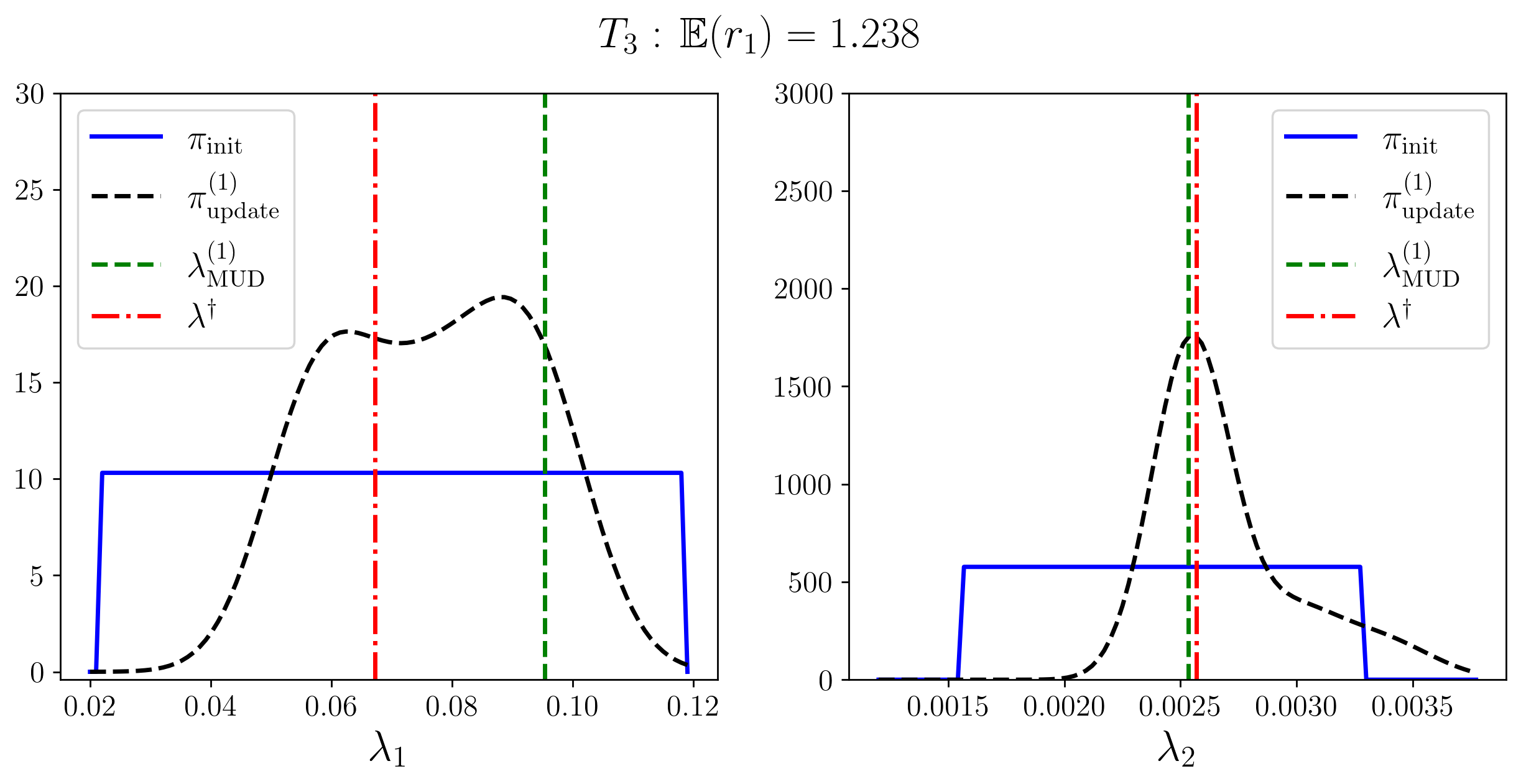}
   \includegraphics[width=0.80\linewidth]{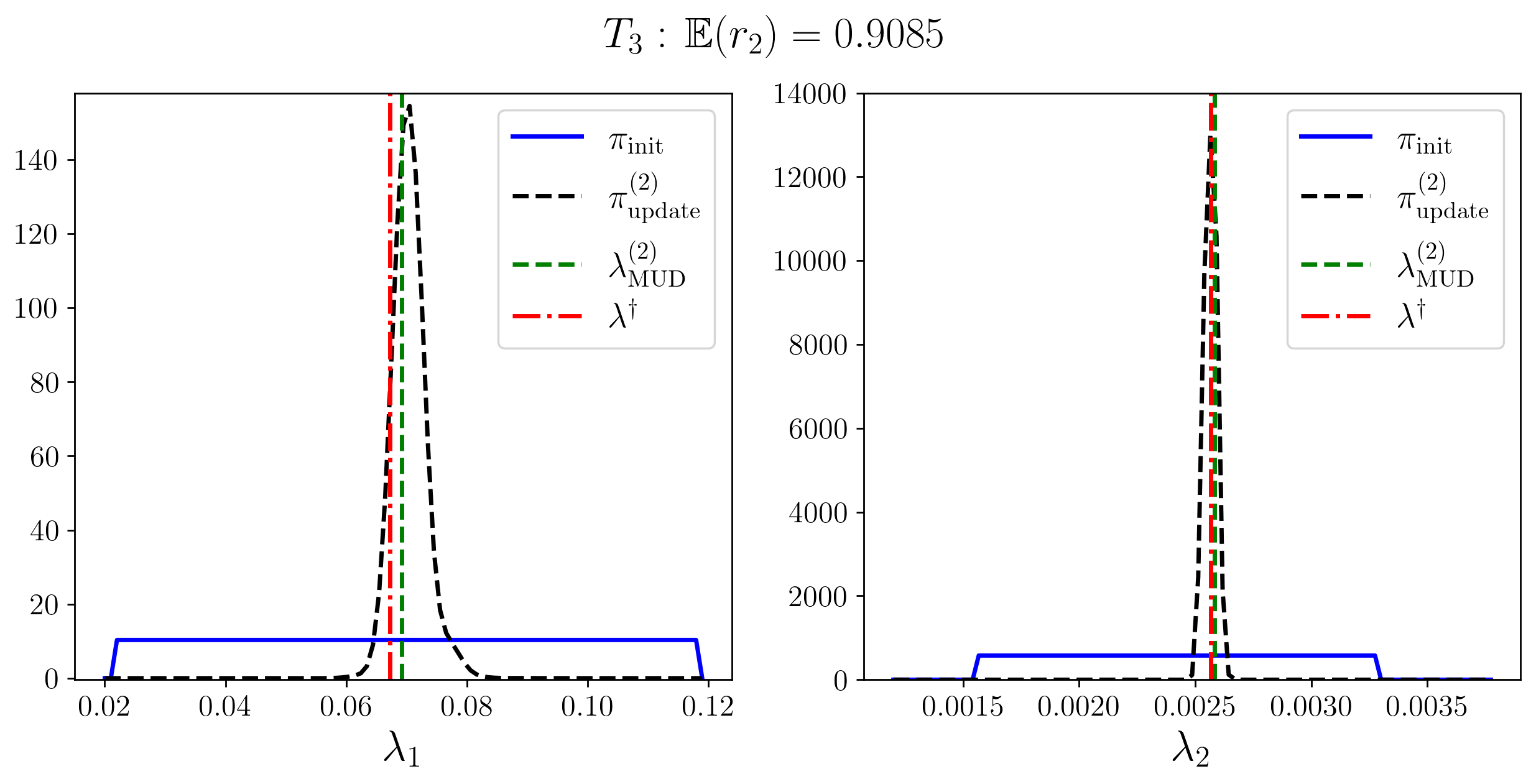}
   \caption{Updated distributions (dotted black line) and mud estimates (dotted green line) using a one component (top) vs two component (bottom) $Q_{PCA}$ map for $T_3$ (see Figure~\ref{fig:adcirc-res}) with $N = 959$ data points from January 07 00:00:00 to January 09 00:00:00.}
   \label{fig:adcirc-t3}
\end{figure}

Three different time windows of data are used to construct $Q_{PCA}$.
In each example, the diagnostic $\mathbb{E}(r)$ is used to compare using one ($\mathbb{E}(r_1)$) vs. two ($\mathbb{E}(r_2)$) principal components.
First in $T_1$ (Figure ~\ref{fig:adcirc-t1}), we choose a time window with low winds and little variation in time-series values consisting of 119 data points.
Note the $Q_{PCA}$ map with one component (top row) leads to much better $\mathbb{E}(r)$ ($0.9879$ vs $0.4204$).
This makes sense because during a period of low winds the cut-off parameter $\param_2$ will not come affect the resulting dynamics of the system for most samples.
Subsequently, attempting to update both parameters via a two-component QoI map leads to a poor $\mathbb{E}(r)$ indicating a violation of the predictability assumption and that the updated distribution should not be trusted.
We see similar results in $T_2$ (Figure ~\ref{fig:adcirc-t2}), a time window of 53 data points characterized by high winds and large variations in time-series values, where a one-component $Q_{PCA}$ map also performs better in terms of the diagnostic $\mathbb{E}(r)$ ($1.039$ vs $4.42$).
However in $T_2$, we observe that the data are now sensitive to $\param_2$ and not $\param_1$ due to the fact the cut-off parameter dominates impacts the dynamics in a period of high winds. 
Finally, in $T_3$ (Figure ~\ref{fig:adcirc-t3}), a larger time window with both high and low winds is used consisting of 959 data points.
Here, enough data are collected that are sensitive to perturbations in both parameters so that the two-component $Q_{PCA}$ map produces a reasonable diagonstic value $\mathbb{E}(r_2)=0.9085$.
We note that the update to the initial distributions for both parameters is significant in this case, with the two-component $Q_{PCA}$ leading to accurate MUD estimates of the true parameter values.

To summarize, these results on different time windows illustrate well the benefits of using the $Q_{PCA}$ and MUD point estimation  algorithm for parameter estimation problems based on temporal data.
It is critical to monitor the diagnostic $\mathbb{E}(r)$ since it gives us a specific metric to determine the quality of reconstructed distributions and potential violations of assumptions.



\section{Conclusions and Future Work}\label{sec:conclusions}

A new approach to estimating and quantifying uncertainties in estimates of parameter identification problems is presented within a data-consistent framework for solving stochastic inverse problems.
This approach identifies the parameter that maximizes the updated density solving data-consistent inverse problems.
This parameter, referred to as the MUD point, is compared to the maximum a posteriori (MAP) and least squares estimates obtained by solving other formulations of the parameter identification problem.
Under the standard assumptions of linear maps with Gaussian distributions, it is demonstrated that the MUD point maintains the predictive accuracy of a least squares estimate with the flexibility of incorporating prior/initial beliefs that make the MAP estimate popular in the uncertainty quantification literature.
The theory of existence and uniqueness of the MUD point is also proven under these assumptions.

We then demonstrate how to formulate and apply a data-constructed QoI map for estimating and quantifying uncertainties in the MUD point.
Specifically, the definition of a weighted mean error (WME) QoI map is utilized that can incorporate an arbitrary stream of data associated with repeated measurements.
The WME map construction is motivated based on an application of the Fisher-Newman factorization theorem to the joint data-likelihood function commonly used in a Bayesian approach.
The particular form of the WME map used in this work admits a fixed observed density as a function of observed data that nonetheless possesses several useful characteristics. 
First, for linear measurement operators, it is established that once a threshold on the number of observed data is reached, the existence and uniqueness of the MUD point is guaranteed.
Second, the eigenvalues of the updated covariance related to the data-informed directions decrease at rates inversely proportional to the number of data obtained for each component of the WME map.

We then illustrate how an alternative data-constructed QoI map utilizing a PCA can be substituted for the WME map when potentially non-repeated spatial-temporal measurement data are available.
This alternative QoI map is applied to two differential equations utilizing noisy spatial and temporal data including an example where wind drag parameters are estimated within the ADCIRC model for storm surge based on water elevation data over time. 

Future directions include issues of optimal experimental design within this framework as well as the sequential estimation of the MUD point as data are obtained sequentially in space or time.
In particular, for an online learning situation where optimal parameters are to be determined as data ``stream into'' the system, we will consider the construction of the $Q_{PCA}$ maps using increasing time windows of data, with the diagnostic $\mathbb{E}(r)$ serving as a measure as to when enough data have been accumulated to invert on one or more model parameter values.
Re-sampling can then be done from these updated distributions and the process repeated to increase the resolution of the $\mudpt$ estimate as needed.
We will also consider non-density based approaches in future work based on optimization formulations of empirical distribution functions where too few samples are available for reliable density estimation in the space defined by the data-derived QoI. 
Preliminary work in this direction is promising (e.g., see \cite{bergstrom2020cdfs}).
There is also interest in applying this MUD and data-derived QoI framework to higher-dimensional parameter spaces modeling microstructures in material properties as was previously analyzed in \cite{tran2021solving} the classic data-consistent inverse framework.


\section* {ACKNOWLEDGMENTS}\label{sec:acknowledge}
T.~Butler's, M.~Pilosov's, and T.~Yu Yen's work is supported in part by the National Science Foundation under Grant No.~DMS-1818941.
T.~Butler's work is also partially supported by the National Science Foundation under Grant No.~DMS-2208460.
C.~Dawson's and C.~del-Castillo-Negrete's work is supported in part by the National Science Foundation No.~DMS-1818847 and No.~DMS-2208461.
The opinions, findings, and conclusions or recommendations expressed in this manuscript are those of the authors and do not necessarily reflect the views of the National Science Foundation, nor Sandia National Laboratories (SAND2022-16991 O).

\appendix

\section{Derivation of $\lambda^{MUD}$ for Linear Gaussian Case}\label{app:mud_deriv}

We begin by an alternative representation of $J(\param)$.
First, define
\begin{equation}\label{eq:eff_reg}
	R := \nolinebreak \initialCov^{-1} - \nolinebreak A^\top \predictedCov^{-1} A.
\end{equation}
Using this $R$, rewrite $J(\param)$ as
\begin{equation}\label{eq:dci-objective-alt}
J(\param):= \norm{\observedMean - Q(\param)}_{\observedCov^{-1}}^2 + \norm{\param - \param_0}_{R}^2.
\end{equation}
In this form, we identify $R$ as the {\em effective regularization} in $J(\param)$ due to the formulation in the data-consistent framework.

Observe that if $m=p$, then, by the assumption that $A$ is full-rank, $A$ is invertible.
In this case, $R$ is the $p\times p$ zero matrix and~\eqref{eq:dci-objective-alt} reduces to the data-discrepancy term so that the MUD point is recognizable as the least squares solution, i.e., the point that minimizes the data-discrepancy term.
Moreover, in this case we can immediately identify that $\mudpt = A^{-1}(\observedMean-b)$.
This is also evident from the perspective of the densities.
Specifically, in this case, $\updated$ is defined by applying a change of variables formula to $\observed$.

Suppose instead that $m<p$ so that the inverse-problem is under-determined.
In this case, we observe that constructing $R$ only requires specification of the initial/prior density and the QoI map, i.e., $R$ may be defined prior to any collection of data on the QoI.
Subsequently, we can interpret $J(\param)$ as coming from a modified Bayesian inverse problem with a prior defined by a $N(\param_0,\Sigma_R)$ distribution where $\Sigma_R=R^{-1}$.
In other words, the MUD and MAP points can both be interpreted as solutions to different Bayesian inverse problems.

However, $\Sigma_R$ is in fact a degenerative covariance, i.e., $R$ is not technically invertible.
This implies that $\Sigma_R$ cannot be directly substituted in for $\initialCov$ in~\eqref{eq:map_cov_analytical} to define a closed form expression for $\updatedCov$.
We therefore first substitute  $\Sigma_\text{post}$ and $\initialCov^{-1}$ in~\eqref{eq:map_cov_analytical} with $\updatedCov$  and $R$, respectively, to get
\begin{equation}
	\updatedCov := \left(A^\top \observedCov^{-1} A + R\right)^{-1}.
\end{equation}
Since $R$ is not invertible, Woodbury's identity cannot be applied (yet).
Using~\eqref{eq:eff_reg}, we can form
\begin{equation}
	\updatedCov = \left(A^\top \observedCov^{-1} A +  \initialCov^{-1} - \nolinebreak A^\top \predictedCov^{-1} A\right)^{-1},
\end{equation}
which is re-arranged as
\begin{equation}
	\updatedCov = \left(A^\top \left[\observedCov^{-1} - \predictedCov^{-1}\right]A + \initialCov^{-1}\right)^{-1}.
\end{equation}
Recall from Section~\ref{subsec:Motivation} that the predictability assumption in this case is that the smallest eigenvalue of $\predictedCov$ is larger than the largest eigenvalue of $\observedCov$.
The roles are reversed when we consider the inverses of these matrices.
Subsequently, $\observedCov^{-1}-\predictedCov^{-1}$ is a symmetric positive definite matrix and thus invertible.
Applying the Woodbury identity yields
\begin{equation}\label{eq:updated_cov_almost}
	\updatedCov = \initialCov - \initialCov A^\top\left( \left[\observedCov^{-1} - \predictedCov^{-1}\right]^{-1} + \predictedCov\right)^{-1} A\initialCov.
\end{equation}
Applying Hua's identity and simplifying gives
\begin{equation}\label{eq:Hua}
	\left( \left[\observedCov^{-1} - \predictedCov^{-1}\right]^{-1} + \predictedCov\right)^{-1} = \predictedCov^{-1}\left[\predictedCov - \observedCov\right]\predictedCov^{-1}.
\end{equation}

Substituting~\eqref{eq:Hua} into~\eqref{eq:updated_cov_almost} gives
\begin{equation}\label{eq:updatedCov_final}
	\updatedCov = \initialCov - \initialCov A^\top \predictedCov^{-1}\left[\predictedCov-\observedCov\right]\predictedCov^{-1}A\initialCov.
\end{equation}
We can now modify the expression for the MAP point given in~\eqref{eq:map-point-analytical} by substituting $\updatedCov$ for $\Sigma_\text{post}$ to write the MUD point that minimizes $J$ as
\begin{equation}\label{eq:mud-point-analytical-alt}
\param^{\text{MUD}} = \param_0 + \updatedCov A^\top \observedCov^{-1} (\observedMean - b - A\param_0).
\end{equation}
Substituting~\eqref{eq:updatedCov_final} into~\eqref{eq:mud-point-analytical-alt} and simplifying, we have
\begin{equation}\label{eq:mud-point-analytical-appendix}
	\mudpt = \param_0 + \initialCov A^\top \predictedCov^{-1}(\observedMean - b - A\param_0).
\end{equation}

From a practical perspective, ~\eqref{eq:mud-point-analytical-appendix} is the preferred form for calculating the MUD point numerically given its reduced complexity in terms of the number of FLOPS required, and it is the default method used in the software implementation (see below).
One would opt in for using ~\eqref{eq:mud-point-analytical-alt} only if the updated covariance $\updatedCov$ is required.
This option is available in the software implementation (See \ref{app:A}) by setting the \verb|solve(method=`mud_alt')| option in the \verb|LinearGaussianProblem| class.
The updated covariance $\updatedCov$ is computed using the \\
\verb|LinearGaussianProblem.updated_cov()|function which computes ~\eqref{eq:updatedCov_final}.

\section{Software Contributions}\label{app:A}

The work presented here is available on GitHub at \verb|github.com/mathematicalmichael/mud.git| as a complete Python package compliant with PEP 517 and 518 and published to the PyPi Python Package Registry under the name \verb|mud|.
Convenient python classes encompassing the core mathematical objects and analytical expressions are provided in the \verb|mud.base|, including \verb|LinearGaussianProblem|, for the analytical linear solutions, \verb|DensityProblem|, for the density base MUD point estimation problems, and \verb|SpatioTemporalProblem|, for aggregating spatio-temporal data in data-constructed QoI maps.

Running \verb|pip install mud[examples]| will install the mud package and its dependencies, as well as a convenient Command Line Interface (CLI) to run the main examples presented.
This means that upon successful installation, one can run \verb|mud examples --help| from a command line to explore the examples available and options for each.
The entrypoint \verb|mud examples mud-paper| will run all the examples as shown in this paper (except for the ADCIRC example in section 6.2, see below).
The package version used in this paper is \verb|mud==0.1|, which is compatible with Python 3.7+.

The datasets for the ADCIRC example in section 7.2 is hosted on DesignSafe, a comprehensive cyber-infrastructure that is part of the NSF-funded Natural Hazard Engineering Research Infrastructure (NHERI) and provides cloud-based tools to manage, analyze, understand, and publish critical data and research related to impacts of natural hazards \cite{rathje2017designsafe}.
The published project directory \cite{carlos2022parameter} includes a static version of the ADCIRC data presented, along with Jupyter notebooks to recreate the data-set itself using DesignSafe HPC resources to run ensembles of ADCIRC simulations.


\bibliographystyle{elsarticle-num-names}
\bibliography{ReferencesBib.bib}

\end{document}